\documentclass[reqno,twoside,12pt,a4paper]{amsart}
\usepackage{layout}
\usepackage{latexsym}
\usepackage{amsmath}
\usepackage{amssymb}
\usepackage{cases}
\usepackage{ascmac}
\usepackage{mathrsfs}
\usepackage{amsthm}
\usepackage{tikz}
\usetikzlibrary{intersections, calc, arrows}
\usepackage{color}
\def\pier #1{{\color{blue} #1}}
\def\takeshi #1{{\color{magenta}#1}}
\def\luca #1{{\color{green}#1}}

\def\pier #1{#1}
\def\takeshi #1{#1}
\def\luca #1{#1}

\numberwithin{equation}{section}

\newtheorem{theorem}{Theorem}[section]

\newtheorem{proposition}[theorem]{Proposition}

\newtheorem{remark}[theorem]{Remark}

\allowdisplaybreaks[4]

\title[Cahn--Hilliard system with dynamic 
boundary conditions]
{A Cahn--Hilliard system with forward-backward dynamic 
boundary condition\\ and non-smooth potentials}

\author[P.\ Colli]{Pierluigi Colli}
\address{Pierluigi Colli: Dipartimento di Matematica, Universit\`a degli Studi di Pavia,
Via Ferrata~1, 27100 Pavia, Italy}
\email{pierluigi.colli@unipv.it}

\author[T.\ Fukao]{Takeshi Fukao}
\address{Takeshi Fukao: Department of Mathematics, Faculty of Education, 
Kyoto University of Education, 
1~Fujinomori, Fukakusa, Fushimi-ku, Kyoto~612-8522 Japan}
\email{fukao@kyokyo-u.ac.jp}

\author[L.\ Scarpa]{Luca Scarpa}
\address{Luca Scarpa: Department of Mathematics, Politecnico di Milano, 
Via E.~Bonardi 9, 20133 Milano, Italy}
\email{luca.scarpa@polimi.it}
\urladdr{https://sites.google.com/view/lucascarpa}

\dedicatory{}
\subjclass[2000]{}
\begin{document}

\thispagestyle{empty}
\maketitle

\begin{abstract}
A system with equation and dynamic boundary condition of Cahn--Hilliard type is considered. 
This system comes from a derivation performed in Liu--Wu \takeshi{
(Arch.\ Ration.\ Mech.\ Anal.\ {\bf 233} (2019), 167--247) via an energetic variational approach.} 
Actually, the related problem can be seen as a transmission problem 
for the phase variable in the bulk and the corresponding variable on the boundary. 
The asymptotic behavior as the coefficient of the surface diffusion acting on the 
boundary phase variable goes to $0$ is investigated. By this analysis we obtain 
a forward-backward dynamic boundary condition at the limit. We can deal with 
a general class of potentials having a double-well structure, including the 
non-smooth double-obstacle potential. We illustrate that the limit problem 
is well-posed by also proving a continuous dependence estimate. 
Moreover, in the case when the two graphs, in the bulk and on the boundary, 
exhibit the same growth, we show that
the solution of the limit problem is more regular and we prove 
an error estimate for a suitable order of the diffusion parameter.
\medskip

\medskip

\end{abstract}

\section{Introduction}
\label{intro}
\setcounter{equation}{0} 

Let $T>0$ be some finite time and let $\Omega \subset \mathbb{R}^d$ ($d=2,3$) be a bounded smooth domain. Consider the heat equation: for a given initial data $u_0:=u_0(x)$, $x\in \Omega$, and heat source $f:=f(t,x)$, find $u:=u(t,x)$, $(t,x) \in Q:= (0,T) \times \Omega$, satisfying 
\begin{equation}
\label{heat1}
	\partial_t u -\Delta u = f \quad \hbox{in~}Q, \quad u(0)=u_0 \quad \hbox{in~}\Omega,
\end{equation}
besides some suitable boundary condition. If instead the sign in front of the {L}aplace term $\Delta u$ appearing in the heat equation is positive, that is, 
\begin{equation}
\label{heat2}
	\partial_t u +\Delta u = f\quad \hbox{in~}Q, \quad u(0)=u_0 \quad \hbox{in~}\Omega,
\end{equation}
the resultant is known to be an ill-posed problem. 
Indeed, \luca{\eqref{heat2} is backward-in-time} and can be interpreted as a determination problem of the history of heat diffusion as follows: 
by the change of variable $U(t):=u(T-t)$, $t\in (0,T)$, we obtain
\begin{equation}
\label{heat3}	\partial_t U - \Delta U =-f\quad \hbox{in~}Q, \quad U(T)=u_0\quad \hbox{in~}\Omega,
\end{equation}
where the initial condition is changed as a terminal condition at time $T$.
From the general theory of partial differential equations, 
it is known that the forward heat equation \eqref{heat1} has the special property of the smoothing effect. 
More precisely, you can gain the smoothness of the solution at any short time even if the initial datum is not so smooth. 
Therefore, this consideration suggests us that some small noise in the terminal data may come from pathological deviations on intermediate states for the backward heat equation \eqref{heat3}. 
In this sense, the continuous dependence is a delicate problem and  we can say that the backward heat equation is ill posed, in general. 
The issue of the existence of solutions is also delicate. In order to discuss it, one needs some additional settings (see, e.g., \cite{Mir61}).

About this class of problems, let us raise the question: what can happen when 
the backward problem is set on the boundary as a dynamic boundary condition?

In this paper, we are concerned with a (possible) backward heat equation on the boundary $\Gamma :=\partial \Omega$ of some smooth bounded domain $\Omega \subset \mathbb{R}^d$ 
($d=2,3$); namely, we address a backward equation as a dynamic boundary condition of a problem which consists in finding $v : \Sigma \to \mathbb{R}$ that satisfy
\begin{align*}
	\partial _t v  + \Delta_\Gamma v = \takeshi{Gu}
	& \quad {\rm on~} \Sigma:=(0,T) \times \Gamma, \\
	v(0) = v_0 & \quad {\rm on~} \Gamma,
\end{align*}
where $\partial_t$ and $\Delta_\Gamma$ stand for the time derivative and 
{L}aplace--{B}eltrami operator (see, e.g., \cite{Gri09}), respectively.  
Moreover, $v_0:\Gamma \to \mathbb{R}$ is prescribed. 
\luca{The backward nature of the boundary problem is due to the fact} that 
the sign of the {L}aplace--{B}eltrami  term appearing in the 
dynamic boundary condition is positive. 
The detail about the right-hand side \takeshi{$Gu$} is given later: indeed, the variable
$u:Q:=(0,T) \times \Omega \to \mathbb{R}$ is also unknown and runs in the bulk, being 
related by a transmission condition to the unknown $v : \Sigma \to \mathbb{R}$
on the boundary.

\luca{In order to give rigorous sense to the backward dynamics on the boundary,
first we artificially provide the problem with a suitable equation in the bulk
with a fourth-order boundary condition in such a way that 
the respective bulk-boundary problem is well-posed. 
Then, by performing a vanishing diffusion on the boundary,
in particular we recover the second-order backward heat equation on the boundary.}
The equation considered in the bulk is of {C}ahn--{H}illiard type (see \cite{CH58}),
that refers to a celebrated model describing the spinodal decomposition 
in a simple framework of fourth-order partial differential equations. 
Some historical and mathematical description of {C}ahn--{H}illiard systems can be found in the papers \cite{Cah11, EZ86, KNP95, Mir17, Mir19}, to mention only a few. \luca{On the boundary,}
we consider the following dynamic condition of {C}ahn--{H}illiard type 
(see, e.g., \cite{CFW20, GK, LW19}): for $\delta \in (0,1]$ we look for 
$v : \Sigma \to \mathbb{R}$ fulfilling 
\begin{align}
  \partial _t v  - \Delta_\Gamma w = 0
  & \quad {\rm on~} \Sigma, \label{lw5}\\ 
  w= - \delta \Delta_\Gamma v+\beta_\Gamma (v)+\pi_\Gamma(v)-g+\partial_{\boldsymbol{\nu}} u
  & \quad {\rm on~} \Sigma, \label{lw6}\\
  v(0) = v_0 & \quad {\rm on~} \Gamma, 
  \label{lw8}
\end{align}
where 
$\beta_\Gamma$ is a monotone function (it may be also a graph), 
$\pi_\Gamma$ is an anti-monotone Lipschitz continuous function, $\partial _{\boldsymbol{\nu}}$ stands for the 
normal derivative, $g: \Sigma \to \mathbb{R}$ is a given datum. 
In the last term of \eqref{lw6} the normal derivative of another unknown function 
$u: Q \to \mathbb{R}$ appears, and correspondingly $u$ has to satisfy 
\begin{align}
  \partial _t u - \Delta \mu =0 
  & \quad {\rm in ~ } Q, \label{lw1}\\
  \mu = -\Delta u +\beta(u)+\pi(u)- f 
  & \quad {\rm in~} Q, \label{lw2}\\
  \partial _{\boldsymbol{\nu}} \mu =0 
  & \quad {\rm on~} \Sigma, \label{lw3}\\ 
  u_{|_\Gamma} = v
  & \quad {\rm on~} \Sigma, \label{lw4}\\ 
  u(0) = u_0 & \quad {\rm in~} \Omega,
  \label{lw7}
\end{align}
where the symbol $\Delta$ stands for the {L}aplacian, $u_{|_\Gamma}$ represents the trace of $u$ on $\Gamma$, $\beta$ and $\pi $
play the same role in the bulk as $\beta_\Gamma$ and $\pi_\Gamma$ on the boundary,
$f: Q \to \mathbb{R}$ is another datum. Of course, in \eqref{lw5}--\eqref{lw7}
two auxiliary variables $w: \Sigma \to \mathbb{R}$ and $\mu: Q \to \mathbb{R}$, 
which have the physical meaning of chemical potentials, are also outlined.
 
Here, we intentionally construct the system from the equations 
on the boundary with side conditions on the bulk.  
This implies that the system presents the main equations
on the boundary with the equations in the bulk interpreted 
as auxiliary conditions (same procedure as, e.g., in \cite{CF20, FIK12, FIKL20} and references therein). 
\takeshi{Note that if we simply take $\beta _\Gamma(r)=0$, 
$\pi_\Gamma(r)=-r$ for $r \in \mathbb{R}$, 
and let $\delta \to 0$ in \eqref{lw5}--\eqref{lw6}, 
then the target equation on the boundary reads}
\begin{equation}
	\partial _t v  + \Delta_\Gamma v = \takeshi{Gu}:=\Delta_\Gamma (\partial _{\boldsymbol{\nu}} u -g)
	\quad {\rm on~} \Sigma \label{back}
\end{equation}
and actually \luca{makes sense as a} backward equation. 
On the other hand, the complementary system \eqref{lw1}--\eqref{lw7} 
is ready to help in order to gain solvability of the full problem despite the backward equation on the boundary.
 
The main topic of this paper is related to the rigorous discussion of the limiting 
procedure as $\delta \to 0$ for the complete system~\eqref{lw5}--\eqref{lw7}
and the novelty is the treatment of wider \takeshi{classes for $\beta $ and $\beta_\Gamma$.
Indeed, we can} postulate that $\beta$ and $\beta_\Gamma$ are maximal monotone graphs, that may be  
multivalued, with suitable growth properties. In this respect, the equations \eqref{lw6} and \eqref{lw2} should be rewritten for suitable selections $\eta$ of $\beta_\Gamma (v)$ and $\xi $ of $\beta (u)$, respectively.  In fact, in our approach $\beta$ and $\beta_\Gamma$ are the 
subdiffentials of proper convex lower semicontinuous functions $\widehat{\beta}, \,
\widehat{\beta}_\Gamma:\mathbb{R} \to [0,+\infty]$ such that  
$\widehat{\beta}(0)=\widehat{\beta}_\Gamma(0)=0$, and the growth of $\beta$ is dominated
by the one of $\beta_\Gamma$, in the sense of assumption {\rm (A1)} below with condition  
\eqref{dom2}. In this framework, we can prove that the solution to~\eqref{lw5}--\eqref{lw7}, whose determination is ensured by the results in \cite{CFW20}, suitably converges as $\delta \to 0$ to the solution of the limit problem in which \eqref{lw6} is replaced by the analogous condition with $\delta=0$. Actually, it occurs that 
in the limiting process the solution of the problem with $\delta \in (0,1]$ looses some regularity at the limit, and the limit boundary equation $  w = \partial_{\boldsymbol{\nu}} u  - g +\beta_\Gamma (v) +\pi_\Gamma(v)$
has to be properly interpreted in the sense of a subdifferential 
inclusion in dual spaces. However, the limit problem turns out to exhibit a well-posedness property
since the continuous dependence of the solution with respect to the initial data and the source terms $f$ and $g$ can be proved. In addition to these results, in the special situation when the two graphs $\beta$ and $\beta_\Gamma$ have a comparable growth
(\takeshi{cf.\ assumption} \eqref{ip_extra} later on), we show that the solution enjoys more regularity and the limit boundary equation makes sense \luca{also almost everywhere}. 
Moreover, we examine the refined convergence and arrive at an error estimate, for the difference of solutions, of order $\delta^{1/2}$.

Let us now mention some related work. 
Recently the equation and dynamic boundary condition of {C}ahn--{H}illiard type have been 
studied in several papers from various viewpoints. In particular, the {C}ahn--{H}illiard 
system coupled with the dynamic boundary condition of {C}ahn--{H}illiard type as 
\eqref{lw5}--\eqref{lw7} has been introduced and examined by {L}iu--{W}u in \cite{LW19} 
for smooth or singular potentials. Then, it is important to quote the article~\cite{GK} 
where the same problem is treated with a gradient flow approach.
After that, the well-posedness problem for non-smooth potentials has been discussed 
in \cite{CFW20}. Among other contributions for this model, we point out \cite{MW20} 
for the long time behavior and \cite{Met21} for the numerical analysis. 
As a remark, there is a similar system of equation and dynamic boundary condition 
of {C}ahn--{H}illiard type, which has been analysed, earlier than the one in \cite{LW19}, 
by {G}al \cite{Gal06} or {G}oldstein--{M}iranville--{S}chimperna~\cite{GMS11}. 
For this similar model, which however does not postulate a transmission condition like 
\eqref{lw4}, the same authors of this paper investigated the problem with 
forward-backward boundary condition in \cite{CFS21}. 
A sort of intermediate problem between {G}oldstein--{M}iranville--{S}chimperna~\cite{GMS11} 
and {L}iu--{W}u~\cite{LW19} has been considered 
(see, e.g., \cite{BZ21, KLLM20}). About the vanishing diffusion on the dynamic boundary condition, 
the reader may also see the treatments in \cite{CF20b, Sca19} for other \takeshi{{C}ahn--{H}illiard} systems,
\luca{as well \takeshi{as \cite{BCST18, BCST20, CS16} for} vanishing diffusion in the bulk and convergence to 
regularised forward-backward problems.} In the light of vanishing diffusion, 
let us additionally mention the contributions \cite{CF16, Fuk16},
in which the asymptotic limit of a {C}ahn--{H}illiard system converging to a nonlinear diffusion 
equation is considered: the approach of \cite{CF16, Fuk16} consists in taking, for 
$\delta \in (0,1]$, the {C}ahn--{H}illiard system
\begin{align*}
	\partial _t u  - \Delta \mu = 0
	& \quad {\rm in~} Q, \\
	\mu= - \delta \Delta u + \beta (u)+ \delta \pi(u) - f
	& \quad {\rm in~} Q,
\end{align*}
with {N}eumann boundary conditions, 
where  the functions 
$\beta$ and $\delta \pi$ are the monotone and anti-monotone parts of the derivative of a double well potential.
Letting $\delta \to 0$, the target problem is based on the nonlinear diffusion equation 
$\partial _t u  - \Delta (\beta (u) - f) =0$ in $Q$. Similar asymptotic limits have been applied also in other contexts (see, e.g., \takeshi{\cite{Fuk16b, FKY18, GK, KO22, Kur21, KY17, KY17b, Sca20, TST14}}).

We present a brief outline of the paper which is structured as follows.
In Section~\ref{main}, the reader can find the notation and the basic tools for a precise interpretation of the problem, which is clearly stated in terms of variational equations and regularity of solutions. After that, the main theorems are precisely stated.
Section~\ref{uniform} is devoted to the proof of the uniform estimates, independent of the coefficient $\delta$, for the solution to a viscous approximation of the system~\eqref{lw5}--\eqref{lw7}, this viscous approximation having already been used
in~\cite{CFW20}. Finally, in Section~\ref{proofs} the main theorems are finally proved,
with the proofs presented in this order: we start with proving the passage to the limit 
as $\delta \to 0$ on the basis of the uniform estimates; next, we deal with the continuous dependence estimate, of the solution with respect to the data; then, we examine the refined convergence and show the error estimate of order $\delta^{1/2}$ in the case when 
the two graphs exhibit the same growth.

\section{Main theorems}
\label{main}

In this section, we present the main theorems. 
To this aim, we set up the target problem and its fundamental settings. 

\subsection{Notation and useful tools}

Let $T>0$ be a finite time and let $\Omega \subset \mathbb{R}^d$ ($d=2,3$) be a bounded domain 
with smooth boundary $\Gamma :=\partial \Omega$. Moreover, we define the sets
$Q:=(0,T) \times \Omega$ and $\Sigma:=(0,T)\times \Gamma$. 
We use the following notation for the function spaces: $H:=L^2(\Omega)$, 
$V:=H^1(\Omega)$, and $W:=H^2(\Omega)$.  
Norms and inner products will be denoted by $|\cdot |_X$ and $(\cdot,\cdot)_X$, respectively, 
where $X$ is the corresponding {B}anach or {H}ilbert space.
Analogously, let $H_\Gamma:=L^2(\Gamma)$, $V_\Gamma:=H^1(\Gamma)$, 
$W_\Gamma:=H^2(\Gamma)$, and set $Z_\Gamma:=H^{1/2}(\Gamma)$ as well.
Next, we define the bilinear forms $a:V \times V \to \mathbb{R}$ and $a_\Gamma: V_\Gamma \times V_\Gamma \to \mathbb{R}$ by 
\begin{align*}
	a(z,\tilde{z}) & := \int_\Omega \nabla z \cdot \nabla \tilde{z} \,dx  
	\quad {\rm for~} z, \tilde{z} \in V, \\
	a_\Gamma (z_\Gamma,\tilde{z}_\Gamma) & := \int_\Gamma \nabla _\Gamma z_\Gamma \cdot \nabla _\Gamma \tilde{z}_\Gamma \, d\Gamma
	\quad {\rm for~} z_\Gamma, \tilde{z}_\Gamma \in V_\Gamma,
\end{align*}
where the symbol $\nabla _\Gamma$ stands for the surface gradient. 
Moreover, 
we define two functions $m : V^* \to \mathbb{R}$ and $m_\Gamma: V_\Gamma^* \to \mathbb{R}$ by 
\begin{align*}
	m(z^*) & := \frac{1}{|\Omega|} \langle z^*,1 \rangle_{V^*,V} 
	\quad {\rm for~} z^* \in V^*, \\
	m_\Gamma(z_\Gamma^*) & := \frac{1}{|\Gamma|} \langle z_\Gamma^*,1 \rangle_{V_\Gamma^*,V_\Gamma} 
	\quad {\rm for~} z_\Gamma^* \in V_\Gamma^*, 
\end{align*}
where the symbol $X^*$ stands for the dual spaces of the corresponding {B}anach space $X$,  
\takeshi{$|\Omega|:=\int_\Omega 1\, dx$, and 
$|\Gamma|:=\int_\Gamma 1\, d\Gamma$.} If 
$z^* \in H$, then $m(z^*)$ is the mean value of $z^*$. Analogously, $m_\Gamma(z_\Gamma^*)$ has the same meaning for $z_\Gamma^* \in H_\Gamma$.
Using them, we define $H_0:=H \cap {\rm ker} (m)=\{ z \in H : m (z)=0 \}$, 
$H_{\Gamma,0}:=H_\Gamma \cap {\rm ker} (m_\Gamma)$, 
$V_0:=V \cap H_0$, and $V_{\Gamma,0}:=V_\Gamma \cap H_{\Gamma,0}$ with 
the following inner products 
\begin{align*}
	(z,\tilde{z})_{H_0} & :=(z,\tilde{z})_{H} 
	\quad {\rm for~} z, \tilde{z} \in H_0, \\
	(z,\tilde{z})_{V_0} & :=a(z,\tilde{z})
	\quad {\rm for~} z, \tilde{z} \in V_0, \\
	(z_\Gamma,\tilde{z}_\Gamma)_{H_{\Gamma,0} } & :=(z_\Gamma,\tilde{z}_\Gamma)_{H_\Gamma} 
	\quad {\rm for~} z_\Gamma, \tilde{z}_\Gamma \in H_{\Gamma,0}, \\
	(z_\Gamma,\tilde{z}_\Gamma)_{V_{\Gamma,0}} & :=a_\Gamma (z_\Gamma ,\tilde{z}_\Gamma)
	\quad {\rm for~} z_\Gamma, \tilde{z}_\Gamma \in V_{\Gamma,0}.
\end{align*}
We point out that, owing to the {P}oincar\'e--{W}irtinger inequality, there exists 
a constant $C_{\rm P}>0$ such that 
\begin{alignat}{2}
	|z|_{V}^2 
	&\le C_{\rm P} \Bigl( \bigl| z -m(z) \bigr|_{V_0}^2 + 
	\bigl| m(z) \bigr|^2 \Bigr) \quad&&{\rm for~all~} z \in V, \label{pier3}\\
	|z|_{V}^2 
	&\le C_{\rm P} |z|_{V_0}^2 \quad&&{\rm for~all~} z \in V_0, 
	\label{pier4}\\ 
	|z_\Gamma|_{V_\Gamma}^2 
	&\le C_{\rm P}  \Bigl( \bigl| z_\Gamma -m_\Gamma(z_\Gamma) \bigr|_{V_{\Gamma,0}}^2 
	+ \bigl| m_\Gamma(z_\Gamma) \bigr|^2 \Bigr) 
	\quad&&{\rm for~all~} z \in V_{\Gamma},
	\label{pier5}\\
	|z_\Gamma|_{V_\Gamma}^2 
	&\le C_{\rm P} |z_\Gamma|_{V_{\Gamma,0}}^2 
	\quad&&{\rm for~all~} z_\Gamma \in V_{\Gamma, 0}. \label{pier6}
	\end{alignat} 
Therefore, we can define the bounded linear operators
$F:V_0 \to V_0^*$ and $F_\Gamma :V_{\Gamma,0} \to V_{\Gamma,0}^*$ as follows:
\begin{alignat*}{2}
	\langle F z,\tilde{z} \rangle_{V_0^*,V_0} & :=a(z,\tilde{z})
	\quad&&{\rm for~} z, \tilde{z} \in V_0, \\
	\langle F_\Gamma z_\Gamma,\tilde{z}_\Gamma \rangle _{V_{\Gamma,0}^*,V_{\Gamma,0}} & :=
	a_\Gamma (z_\Gamma ,\tilde{z}_\Gamma)
	\quad&&{\rm for~} z_\Gamma, \tilde{z}_\Gamma \in V_{\Gamma,0},
\end{alignat*}
and observe that $F$ and $F_\Gamma$ are duality mappings. 
Moreover, $Fz =0$ in $V_0^*$ if and only if $z=0$ in $V_0$, that is, $F$ is invertible. 
Analogously, $F_\Gamma$ is also invertible. 
Therefore, we can define the inner products 
\begin{alignat*}{2}
	(z^*,\tilde{z}^*)_{V_0^*} 
	& :=\langle z^*,  F^{-1} \tilde{z}^* \rangle_{V_0^*,V_0} 
	\quad&&{\rm for~} z^*, \tilde{z}^* \in V_0^*, \\
	(z_\Gamma^*,\tilde{z}_\Gamma^*)_{V_{\Gamma,0}^* } 
	& :=\langle z_\Gamma^*, F_\Gamma^{-1} \tilde{z}_\Gamma^* \rangle_{V_{\Gamma,0}^*,V_{\Gamma,0}} 
	\quad&&{\rm for~} z_\Gamma^*, \tilde{z}_\Gamma^* \in V_{\Gamma,0}^*\pier{,}
\end{alignat*}
\takeshi{which give the related norms}
\begin{alignat*}{2}
	|z^*|_{V_0^*} 
	& =\left\{ \int_\Omega |\nabla F^{-1} z^* | \, dx \right\}^{1/2} 
	\quad&&{\rm for~} z^* \in V_0^*, \\
	| z_\Gamma^* |_{V_{\Gamma,0}^* } 
	& =\left\{ \int_\Gamma | \nabla _\Gamma F_\Gamma ^{-1}z_\Gamma^*| \, d\Gamma \right\}^{1/2}
	\quad&&{\rm for~} z_\Gamma^*  \in V_{\Gamma,0}^*.
\end{alignat*}
Finally, we introduce the following norms in $V^*$ and $V_{\Gamma}^*$,
\begin{alignat}{2}
	|z^*|_{*} 
	& =\left\{  \bigl| z^* -m(z^*) \bigr|_{V_0^*}^2 + \bigl| m(z^*) \bigr|^2 \right\}^{1/2} 
	\quad&&{\rm for~} z^* \in V^*, 
	\label{norm}\\
	| z_\Gamma^* |_{\Gamma,* } 
	& =\left\{  \bigl|  z_\Gamma^*-m_\Gamma(z_\Gamma^*) \bigr|_{V_{\Gamma,0}^*}^2 + \bigl| m_\Gamma(z_\Gamma^*) \bigr|^2 \right\}^{1/2}
	\quad&&{\rm for~} z_\Gamma^*  \in V_{\Gamma}^*, \label{norm'}
\end{alignat}
and observe that they are equivalent to the standard induced norms 
$|\cdot |_{V^*}$ of $V^*$ and $|\cdot |_{V_\Gamma^*}$ of $V_\Gamma^*$, respectively.
Then we obtain the following dense and compact embeddings: 
\begin{gather*}
	V \mathop{\hookrightarrow} \mathop{\hookrightarrow} 
	H \mathop{\hookrightarrow} 
	V^*, 	
	\quad 
	V_0 \mathop{\hookrightarrow} \mathop{\hookrightarrow} 
	H_0 \mathop{\hookrightarrow} 
	V_0^*, 	
    \\	
	 V_{\Gamma} \mathop{\hookrightarrow} \mathop{\hookrightarrow} 
	 H_\Gamma \mathop{\hookrightarrow}  
	 V_{\Gamma}^*, 
	 \quad 
	 Z_{\Gamma} \mathop{\hookrightarrow} \mathop{\hookrightarrow} 
	 H_\Gamma \mathop{\hookrightarrow}  
	 V_{\Gamma}^*, 
	 \quad 
	 V_{\Gamma, 0} \mathop{\hookrightarrow} \mathop{\hookrightarrow} 
	 H_{\Gamma,0} \mathop{\hookrightarrow}  
	 V_{\Gamma,0}^*, 
\end{gather*} 
where 
``$\mathop{\hookrightarrow} \mathop{\hookrightarrow} $'' stands for 
the dense and compact embedding. 
\smallskip

For the reader's convenience, we recall useful tools in functional analysis. 
The first tool is related to the trace theorem 
(see, e.g., \cite[Theorem~2.24]{BG87}, \cite[Chapter~2, Theorem~5.7]{Nec67}), 
which states that
there exist unique continuous linear operators 
$\gamma_0 : V \to Z_\Gamma$ and 
$\gamma_1 : W \to Z_\Gamma$ 
such that 
\begin{gather*}
	\gamma_0 z =z_{|_\Gamma} \quad {\rm for~all~} z \in C^\infty(\overline{\Omega}) \cap V, \\
	\gamma_1 z = \partial_{\boldsymbol{\nu}} z \quad {\rm for~all~} z \in C^\infty(\overline{\Omega}) \cap W.
\end{gather*}
Moreover, there exists a positive constant $C_{\rm tr}$ such that
\begin{equation}
 	|\gamma_0 z |_{Z_\Gamma} \le C_{\rm tr}|z|_{V} \quad {\rm for~all~} z \in V. \label{tr0}
\end{equation}

\subsection{Target problem}

Now we set up our target problem of the forward-backward dynamic boundary 
equation along with the bulk condition of {C}ahn--{H}illiard type and considering
non-smooth potentials. 
Find 
$v$, $w$, $\eta : \Sigma \to \mathbb{R}$ and 
$u$, $\mu$, $\xi : Q \to \mathbb{R}$
satisfying
\begin{align}
  \partial _t v  - \Delta_\Gamma w = 0
  & \quad {\rm a.e.\ on~} \Sigma, \label{t1}\\ 
  w= \partial_{\boldsymbol{\nu}} u+\eta +\pi_\Gamma(v)-g, \quad 
  \eta \in \beta_\Gamma(v)
  & \quad {\rm a.e.\ on~} \Sigma, \label{t2}\\
  \partial _t u - \Delta \mu =0 
  & \quad {\rm a.e.\ in ~ } Q, \label{t3}\\
  \mu = -\Delta u + \xi +\pi(u)- f, \quad \xi \in \beta (u) 
  & \quad {\rm a.e.\ in~} Q, \label{t4}\\
  \partial _{\boldsymbol{\nu}} \mu =0 
  & \quad {\rm a.e.\ on~} \Sigma, \label{t5}\\ 
  u_{|_\Gamma} = v
  & \quad {\rm a.e.\ on~} \Sigma, \label{t6}\\ 
  v(0) = v_0 & \quad {\rm a.e.\ on~} \Gamma, \label{t8}\\
  u(0) = u_0 & \quad {\rm a.e.\ in~} \Omega, \label{t7}  
\end{align}
where $\beta_\Gamma$ and $\beta$ are maximal monotone graphs on $\mathbb{R} \times \mathbb{R}$, 
$\pi_\Gamma$ and $\pi$ are {L}ipschitz continuous functions, $g:\Sigma \to \mathbb{R}$, $f:Q \to \mathbb{R}$, 
$v_0 : \Gamma \to \mathbb{R}$, and 
$u_0 : \Omega \to \mathbb{R}$ are given functions. Combining \eqref{t1} and \eqref{t2}, 
we find a structure of second order partial differential equation of forward-backward type on the boundary equation. Indeed, in general the sum 
$\beta_\Gamma+\pi _\Gamma$ is not monotonically increasing on the whole domain.
As prototypes, we can choose\smallskip
\renewcommand{\labelitemi}{$\triangleright$}

\begin{itemize}

 \item $\beta_\Gamma(r):=r^3$, $\pi_\Gamma(r):=-r$ for $r \in \mathbb{R}$ (corresponding to the smooth double well potential);\smallskip
 
 \item $\beta_\Gamma(r):=\ln((1+r)/(1-r))$, $\pi_\Gamma(r):=-2cr$ for $r \in (-1,1)$ (derived from the singular potential of logarithmic type, where $c>0$ is a large constant which breaks monotonicity);\smallskip
 
\item $\beta_\Gamma(r):=\partial I_{[-1,1]}(r)$, $\pi_\Gamma(r):=-r$ for $r \in [-1,1]$ (for the non-smooth potential, where the symbol $\partial$ stands for the subdifferential in $\mathbb{R}$);\smallskip
    
 \item $\beta_\Gamma(r):=0$, $\pi_\Gamma (r):=-r$ for $r \in \mathbb{R}$ (for the backward-like heat equation on the boundary).\smallskip
 
\end{itemize}
In our approach, according to previous contributions (cf., e.g., \cite{CC13, CF20b, CFS21, CFW20}),
about $\beta$ we prescribe a condition on the growth, that sets a control 
by the growth of $\beta_\Gamma$, see the later assumption {\rm (A1)} and condition \eqref{dom2}. 
Instead, we can choose any {L}ipschitz continuous function for $\pi$, independent of $\pi_\Gamma$.

\subsection{Main theorems}

We recall an auxiliary {C}ahn--{H}illiard system approaching our target problem: 
for $\delta \in (0,1]$, find $u_\delta$, $\mu_\delta$, $\xi_\delta: Q \to \mathbb{R}$ and 
$v_\delta$, $w_\delta$, $\eta_\delta : \Sigma \to \mathbb{R}$ satisfying
\begin{align}
  \partial _t u_\delta - \Delta \mu_\delta =0 
  & \quad {\rm a.e.\ in ~ } Q, \label{lw1a}\\
  \mu_\delta = -\Delta u_\delta + \xi_\delta +\pi(u_\delta)- f, \quad \xi _\delta\in \beta (u_\delta) 
  & \quad {\rm a.e.\ in~} Q, \label{lw2a}\\
  \partial _{\boldsymbol{\nu}} \mu_\delta =0 
  & \quad {\rm a.e.\ on~} \Sigma, \label{lw3a}\\ 
  (u_\delta)_{|_\Gamma} = v_\delta
  & \quad {\rm a.e.\ on~} \Sigma, \label{lw4a}\\ 
  \partial _t v _\delta - \Delta_\Gamma w_\delta= 0
  & \quad {\rm a.e.\ on~} \Sigma, \label{lw5a}\\ 
  w_\delta= \partial_{\boldsymbol{\nu}} u_\delta- \delta \Delta_\Gamma v_\delta+\eta_\delta +\pi_\Gamma(v_\delta)-g, \quad 
  \eta_\delta \in \beta_\Gamma(v_\delta)
  & \quad {\rm a.e.\ on~} \Sigma, \label{lw6a}\\
   u_\delta(0) = u_0 & \quad {\rm a.e.\ in~} \Omega,
  \label{lw7a} \\
  v_\delta(0) = v_0 & \quad {\rm a.e.\ on~} \Gamma.
  \label{lw8a}
\end{align}
This system of equation and dynamic boundary condition of {C}ahn--{H}illiard type 
has been introduced by {L}iu--{W}u in \cite{LW19} and its solvability is 
discussed in the papers \cite{GK, LW19} under some restrictions 
for $\beta$ and $\beta_\Gamma$, while in the case $\delta>0$ the well-posedness issue 
is examined in \cite[Theorems~2.3, 2.4, and 4.1]{CFW20} under our general conditions 
on the graphs $\beta$ and $\beta_\Gamma$ (cf. the assumption {\rm (A1)} below). 
The aim of the present paper is the extension of the results in \cite{CFW20} to the 
limiting situation $\delta =0$. 
\pier{In particular, in our analysis we are able to avoid the geometric conditions of {L}iu--{W}u (cf.~\cite[Theorem 3.2]{LW19}).}
\smallskip

In this paper, we assume: 
\begin{enumerate}
\item[(A1)]  $\beta$ and $\beta_\Gamma$ are maximal monotone graphs on 
$\mathbb{R} \times \mathbb{R}$, and there exist 
proper, lower semicontinuous, and convex functions 
$\widehat{\beta}, \widehat{\beta}_\Gamma:\mathbb{R} \to [0,+\infty]$ satisfying 
$\widehat{\beta}(0)=\widehat{\beta}_\Gamma(0)=0$ 
 such that 
\begin{equation*}
	\beta =\partial \widehat{\beta}, 
	\quad \beta_\Gamma=\partial \widehat{\beta}_\Gamma.
\end{equation*}
Moreover, we assume that $D(\beta_\Gamma)\subset D(\beta)$ and there exists 
positive constants $\varrho_1, c_1>0$ such that 
\begin{equation}\label{dom2}
  \takeshi{\bigl|}\beta^\circ(r)\takeshi{\bigr|} 
  \leq \varrho_1 \takeshi{\bigl|}\beta_\Gamma^\circ(r)\takeshi{\bigr|} + c_1
  \quad \takeshi{\rm for~all~} \,
  r\in D(\beta_\Gamma).
\end{equation}
\item[(A2)] $\pi, \pi_\Gamma : \mathbb{R} \to \mathbb{R}$ are {L}ipschitz continuous, 
with their constants $L$ and $L_\Gamma$, respectively. 
Moreover, we set $\widehat\pi(\rho):=\int_0^\rho\pi(r)\,dr$ and 
$\widehat\pi_\Gamma(\rho):=\int_0^\rho\pi_\Gamma(r)\,dr$, $\rho\in\mathbb{R}$;   
\item[(A3)] $f \in L^2(0,T;V)$ and $g \in L^2(0,T;V_\Gamma)$; 
\item[(A4)] $u_0 \in V$, $v_0 \in V_\Gamma$ satisfy $\gamma_0 u_0 = v_0$ in $Z_\Gamma$.
Moreover, $u_0\in L^\infty(\Omega)$, \pier{so that $v_0 \in L^\infty(\Gamma)$ as well,} and 
\begin{align*}
\left[\takeshi{\mathop{\rm ess~inf}_{x\in\Omega}} \, u_0(x), \, 
\takeshi{\mathop{\rm ess~sup}_{x\in\Omega}} \, u_0(x)\right]&\subset{\rm int}D(\beta),\\
\left[\takeshi{\mathop{\rm ess~inf}_{x\in\Gamma}} \, v_0(x),\, \takeshi{\mathop{\rm ess~sup}_{x\in\Gamma}} \, v_0(x)\right]&\subset{\rm int}D(\beta_\Gamma).
\end{align*}
Note that this implies that
$\widehat{\beta} (u_0) \in L^1(\Omega)$, $\widehat{\beta}_\Gamma(v_0) \in L^1(\Gamma)$,
$m_{0}:=m(u_0) \in {\rm int} D(\beta)$, and
$m_{\Gamma0}:=m_\Gamma(v_0) \in {\rm int} D(\beta_\Gamma)$.
\end{enumerate}
We notice that in {\rm (A1)} 
the symbol $\beta^{\circ}$ stands for the minimal section defined by 
\begin{equation*}
	\beta^{\circ} (r) :=\bigl\{ r^* \in \beta (r) : |r^*|=\min_{s \in \beta(r)} |s| \bigr\},
\end{equation*}
and same definition holds for $\beta_\Gamma^{\circ}$. 
Of course we can choose $\beta(r)=\beta_\Gamma(r)=0$ for $r \in D(\beta_\Gamma):=\mathbb{R}$. 
\smallskip

Recalling the known result in \cite{CFW20} we obtain the following proposition for $\delta \in (0,1]$.
\begin{proposition}\cite[Theorems~2.3, 2.4]{CFW20}
\label{tak0}
	Under the assumptions {\rm (A1)}--{\rm (A4)}, there exists a sextuplet $(u_\delta,\mu_\delta, \xi_\delta, v_\delta, w_\delta, \eta_\delta)$, 
	where $u_\delta$ and $v_\delta$ are uniquely determined, so that
\begin{gather*}
	u_\delta \in H^1(0,T;V^*) \cap L^\infty (0,T;V) \cap L^2(0,T;W), \\
	\mu_\delta \in L^2(0,T;V), \quad \xi_\delta \in L^2(0,T;H), \\
	v_\delta \in H^1(0,T;V_\Gamma^*) \cap L^\infty (0,T;V_\Gamma) \cap L^2(0,T;W_\Gamma), \\
	w_\delta \in L^2(0,T;V_\Gamma), \quad \eta_\delta \in L^2(0,T;H_\Gamma)
\end{gather*}
and they satisfy 
\begin{gather} 
	\langle \partial_t u_\delta , z \rangle_{V^*,V}
	+ \int_\Omega \nabla \mu _\delta \cdot \nabla z \, dx=0 
	\quad 
	\mbox{for all } z \in V, 
	\mbox{ a.e.\ in } (0,T),
	\label{weak1}
	\\
	\mu_\delta = -\Delta u_\delta + \xi_\delta + \pi (u_\delta) -f, \quad \xi _\delta \in \beta (u_\delta)
	\quad \mbox{a.e.\ in }Q, 
	\label{weak2}
	\\
	(u_\delta)_{|_\Gamma} =v_\delta 
	\quad \mbox{a.e.\ on }\Sigma,
	\label{weak3}
	\\ 
	\langle \partial_t v_\delta, z_\Gamma  \rangle_{V_\Gamma^*,V_\Gamma}
	+ \int_{\Gamma} \nabla_\Gamma w_\delta \cdot \nabla_\Gamma z_\Gamma \, d\Gamma = 0 \quad 
	\mbox{for all } z_\Gamma \in V_\Gamma, 
	\mbox{ a.e.\ in } (0,T),
	\label{weak4}
	\\
	w_\delta = \partial_{\boldsymbol{\nu}} u_\delta- \delta \Delta_\Gamma v_\delta 
	+ \eta_\delta + \pi_\Gamma(v_\delta)-g,  \quad \eta_\delta \in \beta_\Gamma (v_\delta)
	\quad \mbox{a.e.\ on }\Sigma, 
	\label{weak5}
	\\
	u_\delta (0)=u_0
	\quad \mbox{a.e.\ in }\Omega, 
	\label{weak6}
	\\
	v_\delta (0)=v_0
	\quad \mbox{a.e.\ on }\Gamma.
	\label{weak7}
\end{gather}
\end{proposition}

We note that, due to the lack of the regularities of time derivatives, the equations 
(\ref{lw1a}) and (\ref{lw5a}) are replaced by the variational formulations (\ref{weak1}) and (\ref{weak4}), respectively. 
Moreover, the boundary condition (\ref{lw3a}) is hidden in the weak form  (\ref{weak1}). 
Here and hereafter we frequently use the notations $z_{|_\Gamma}$ and $\partial_{\boldsymbol{\nu}} z$ in place of $\gamma_0 z$ and $\gamma_1 z$, respectively. 
\smallskip

Our main theorem is stated here: 

\begin{theorem}\label{tak1}
	Under the assumptions {\rm (A1)}--{\rm (A4)}, 
	there exists at least one sextuplet $(u,\mu, \xi, v, w, \eta)$ fulfilling
\begin{gather*}
	u \in H^1(0,T;V^*) \cap L^\infty (0,T;V), \quad \Delta u \in L^2(0,T;H)\\
	\mu \in L^2(0,T;V), \quad \xi \in L^2(0,T;H), \\
	v \in H^1(0,T;V_\Gamma^*) \cap L^\infty (0,T;Z_\Gamma), \\
	w \in L^2(0,T;V_\Gamma), \quad \eta \in L^2(0,T;Z_\Gamma^*)
\end{gather*}
and satisfying \eqref{t1}--\eqref{t7} in the following sense: 
\begin{gather} 
	\langle \partial_t u , z \rangle_{V^*,V}
	+ \int_\Omega \nabla \mu  \cdot \nabla z \,dx=0 
	\quad 
	\mbox{for all } z \in V, 
	\mbox{ a.e.\ in } (0,T),
	\label{main1}
	\\
	\mu = -\Delta u + \xi + \pi (u) -f, \quad \xi \in \beta (u)
	\quad \mbox{a.e.\ in }Q, 
	\label{main2}
	\\
	 u _{|_\Gamma}=v 
	\quad \mbox{a.e.\ on }\Sigma,
	\label{main3}
	\\ 
	\langle \partial_t v, z_\Gamma  \rangle_{V_\Gamma^*,V_\Gamma}
	+ \int_{\Gamma} \nabla_\Gamma w \cdot \nabla_\Gamma z_\Gamma \,d\Gamma = 0 \quad 
	\mbox{for all } z_\Gamma \in V_\Gamma, 
	\mbox{ a.e.\ in } (0,T),
	\label{main4}
	\\
	(w, z_\Gamma)_{H_\Gamma} = 
	\langle \partial_{\boldsymbol{\nu}} u+\eta, z_\Gamma \rangle_{Z_\Gamma^*,Z_\Gamma}   
	+ \bigl( \pi_\Gamma(v)-g,z_\Gamma \bigr)_{\!H_\Gamma}\quad
	\mbox{for all } z_\Gamma \in Z_\Gamma,
	\mbox{ a.e.\ in } (0,T),
	\label{main5}
	\\
	\langle \eta,z_\Gamma-v \rangle_{Z_\Gamma^*,Z_\Gamma } 
	\le \int_\Gamma \widehat{\beta}_\Gamma(z_\Gamma)\, d\Gamma- 
	\int_\Gamma \widehat{\beta}_\Gamma(v)\, d\Gamma \quad 
	\mbox{for all } z_\Gamma \in Z_\Gamma, 
	\mbox{ a.e.\ in } (0,T),
	\label{main6}
	\\
	u (0)=u_0
	\quad \mbox{a.e.\ in }\Omega, 
	\label{main7}
	\\
	v(0)=v_0
	\quad \mbox{a.e.\ on }\Gamma.
	\label{main8}
\end{gather}
Moreover, $(u,\mu, \xi, v, w, \eta)$ is obtained as limit of the family 
$\{(u_\delta,\mu_\delta, \xi_\delta, v_\delta, w_\delta, \eta_\delta)\}_{0<\delta\leq1}$ 
of the sextuplet solutions given by Proposition~\ref{tak0}, 
in the sense that there is a subsequence $\{\delta_k \}_{k \in \mathbb{N}}$ such that, as $k \to +\infty$, 
\begin{align}
	&u_{\delta_k} \to u \quad \takeshi{\it weakly~star~in~} \,H^1(0,T;V^*) \cap L^\infty (0,T;V),
	\label{cvpier1} \\
	&\Delta u_{\delta_k} \to \Delta u \quad \takeshi{\it weakly~in~} \,L^2 (0,T;H),
	\label{cvpier2} \\
	&\partial_{\boldsymbol{\nu}} u_{\delta_k} \to \partial_{\boldsymbol{\nu}} u 
	\quad \takeshi{\it weakly~in~} \,L^2 (0,T;Z_\Gamma^*),
	\label{cvpier3} \\
	&\mu_{\delta_k} \to \mu  \quad \takeshi{\it weakly~in~} \,L^2 (0,T;V),
	\label{cvpier4} \\
	&\xi_{\delta_k} \to \xi \quad \takeshi{\it weakly~in~} \,L^2 (0,T;H),
	\label{cvpier5} \\
	&v_{\delta_k} \to v \quad \takeshi{\it weakly~star~in~} \,H^1(0,T;V_\Gamma^*) \cap L^\infty (0,T;Z_\Gamma),
	\label{cvpier6} \\
	&\delta_k v_{\delta_k} \to 0 \quad \takeshi{\it strongly~in~} \,L^\infty (0,T;V_\Gamma),
	\label{cvpier9}\\
	&w_{\delta_k} \to w  \quad \takeshi{\it weakly~in~} \,L^2 (0,T;V_\Gamma),
	\label{cvpier7} \\
	&\eta_{\delta_k} \to \eta \quad \takeshi{\it weakly~in~} \,L^2 (0,T;V_\Gamma^*),
	\label{cvpier8} \\
	&(-{\delta_k} \Delta_\Gamma v_{\delta_k} +\eta_{\delta_k}) \to \eta \quad \takeshi{\it weakly~in~} 
	\,L^2 (0,T;Z_\Gamma^*).
	\label{cvpier10}
\end{align}
\end{theorem}

\begin{remark}\rm
  \label{rem1}
  About the inequality~\eqref{main6}, we point out that 
  whenever $\eta \in L^2(0,T; H_\Gamma)$ then
  \eqref{main6} is actually equivalent to the inclusion
  $
  \eta \in \beta_\Gamma(v)$ a.e.~on $ \Sigma$,
  or equivalently
  \[
  \eta \in \partial I_\Sigma(v),
  \]
  where
  \[
  I_\Sigma:L^2(0,T; H_\Gamma)\to[0,+\infty],
  \qquad
  I_\Sigma(z_\Gamma):=
  \begin{cases} \displaystyle 
  \takeshi{\int_\Sigma\widehat\beta_\Gamma(z_\Gamma)\, d\Gamma dt} \quad &\text{if } 
  \widehat\beta_\Gamma(z_\Gamma) \in L^1(\Sigma),\\
  +\infty\quad&\text{otherwise}.
  \end{cases}
  \]  
  On the other hand, if we only have $\eta \in L^2(0,T; Z_\Gamma^*)$, then
  \eqref{main6} means that 
  $
  \eta \in \partial J_\Sigma(v),
  $
  where
  \[
  J_\Sigma:L^2(0,T; Z_\Gamma)\to[0,+\infty],
  \qquad
  J_\Sigma(z_\Gamma):=
  \begin{cases}\displaystyle 
  \takeshi{
  \int_\Sigma\widehat\beta_\Gamma(z_\Gamma) \, d\Gamma dt} \quad&\text{if } 
  \widehat\beta_\Gamma(z_\Gamma) \in L^1(\Sigma),\\
  +\infty\quad&\text{otherwise}.
  \end{cases}
  \]
  Here, the main point is that,
  since we are identifying $H_\Gamma$ to its dual,
  the subdifferential $\partial I_\Sigma$
  is intended as a multivalued operator 
  \[
  \partial I_\Sigma \ \hbox{ from } \ L^2(0,T; H_\Gamma) \ \hbox{ to } \ L^2(0,T; H_\Gamma),
  \]
  while $\partial J_\Sigma$ is seen as an operator, multivalued as well, 
  \[
  \partial J_\Sigma \ \hbox{ from } \ L^2(0,T; Z_\Gamma) \ \hbox{ to } \ L^2(0,T; Z^*_\Gamma).
  \]
For further details we refer to \cite{Bar10, Bre73}.
\end{remark}

\begin{remark}\rm
\label{rem2}
Take, for instance, the case $\beta_\Gamma \equiv 0$, which yields that  $\beta $ should be at most bounded due to {\rm (A1)} and \eqref{dom2}. In this case, it is compulsory to have  $\eta=0$ and,
therefore, by a comparison of term in \eqref{main5} we deduce that $\partial_{\boldsymbol{\nu}} u \in L^2(0,T;Z_\Gamma)$, being in fact
$w=\partial_{\boldsymbol{\nu}} u+ \pi_\Gamma(v)-g $ an element of $ L^2(0,T;V_\Gamma)$. 
Then we interpret the backward equation \eqref{main4} on the boundary as 
\begin{equation*} 
	\langle \partial_t v, z_\Gamma  \rangle_{V_\Gamma^*,V_\Gamma}
	+ \int_{\Gamma} \nabla_\Gamma \takeshi{\bigl(} \partial_{\boldsymbol{\nu}} u
	+ \pi_\Gamma(v) \takeshi{\bigr)}\cdot \nabla_\Gamma z_\Gamma \, d\Gamma = \int_{\Gamma} \nabla_\Gamma g \cdot \nabla_\Gamma z_\Gamma \, d\Gamma	\quad 
	\mbox{for all } z_\Gamma \in V_\Gamma, 
\end{equation*}
a.e.\ in $(0,T)$, where thanks to {\rm (A3)} we can move the term containing $g$ to the right-hand side, but 
we cannot split $\partial_{\boldsymbol{\nu}} u+ \pi_\Gamma(v)\in L^2(0,T;V_\Gamma)$.
\end{remark}

Next theorem is related to the continuous dependence on the given data:
\begin{theorem}
\label{tak2}
For any data $\{ (f^{(i)}, g^{(i)}, u_{0}^{(i)}, v_{0}^{(i)}) \}_{i=1,2}$ 
satisfying {\rm (A3), (A4)} and such that
$m(u_{0}^{(1)})=m(u_{0}^{(2)})$, $m_\Gamma(v_{0}^{(1)})=m_\Gamma(v_{0}^{(2)})$, let 
$(u^{(i)},\mu^{(i)}, \xi^{(i)}, v^{(i)}, w^{(i)}, \eta^{(i)})$ be some respective solutions obtained by 
Theorem~\ref{tak1}. Then there exists a positive constant $C>0$ such that 
\begin{align*}
	& \bigl| u^{(1)}(t) - u^{(2)}(t) \bigr|_*^2 + \bigl| v^{(1)}(t) - v^{(2)}(t) \bigr|_{\Gamma,*}^2
	\nonumber \\
	& \quad {}+ \int_0^t \bigl| u^{(1)}(s) - u^{(2)}(s) \bigr|_{V}^2 \, ds 
	+ \int_0^t \bigl| v^{(1)}(s) - v^{(2)}(s) \bigr|_{Z_\Gamma}^2 \, ds 
	\nonumber \\
	& \le  C \biggl( 
	\bigl| u_0^{(1)} - u_0^{(2)}\bigr|_*^2 + \bigl| v_0^{(1)}- v_0^{(2)} \bigr|_{\Gamma,*}^2
	\nonumber \\
	& \quad {} 
	+ \int_0^t \bigl| f^{(1)}(s) - f^{(2)}(s) \bigr|_{H}^2 \, ds 
	+ \int_0^t \bigl| g^{(1)}(s) - g^{(2)}(s) \bigr|_{H_\Gamma}^2 \, ds \biggr)
\end{align*}
for all $t \in [0,T]$. 
\end{theorem}
Of course, this theorem entails the uniqueness property for $u$ and $v$. If $\beta$ and $\beta_\Gamma$ are single-valued functions, then 
the whole sextuplet $(u,\mu, \xi, v, w, \eta)$ obtained by Theorem~\ref{tak1} is unique as well. 

As a remark, the discussion of the continuous dependence is delicate for backward 
problems in general. In such a problem, under the assumption of the existence of bounded 
solutions, the conditional stability is discussed in some sense in \cite{IY14} (see references therein) and in \cite{SL21} for the \takeshi{{C}ahn--{H}illiard} equation.

The results that follow are inspired by the analogous ones in \cite{CFS21}.

\begin{theorem}
  \label{tak3}
  Under the assumptions {\rm (A1)}--{\rm (A4)}, suppose also that 
  \begin{align}\label{ip_extra}
  &D(\beta)=D(\beta_\Gamma), \quad \hbox{there exists a 
  constant $M\geq1$ such that} 
  \nonumber\\
  &\quad{}\frac1M \takeshi{\bigl|} \beta_\Gamma^\circ(r) \takeshi{\bigr|} 
  - M\leq \takeshi{\bigl|}\beta^\circ(r)\takeshi{\bigr|}  \leq 
  M \takeshi{\bigl( \bigl|} \beta_\Gamma^\circ(r) \takeshi{\bigr|} +1 \takeshi{\bigr)}
  \quad \takeshi{\it for~all~} r\in D(\beta).
  \end{align}
  Then, the limiting sextuplet $(u,\mu, \xi, v, w, \eta)$ 
  obtained in Theorem~\ref{tak1}
  also satisfies
  \begin{align*}
  &u\in L^2 \takeshi{\bigl(} 0,T; H^{3/2}(\Omega) \takeshi{\bigr)}, 
  \ \quad\partial_\nu u \in L^2(0,T; H_\Gamma), \ \quad v\in L^2(0,T; V_\Gamma),\\
  &\eta\in L^2(0,T; H_\Gamma), \quad \ 
  \eta\in\beta_\Gamma(v) \quad\text{a.e.~on } \Sigma.
  \end{align*}
  Moreover, in addition to \eqref{cvpier1}--\eqref{cvpier10}, the following convergences hold, as $k \to +\infty$, 
  \begin{align}
	&\eta_{\delta_k} \to \eta \quad \takeshi{\it weakly~in~} \,L^2(0,T; H_\Gamma),
	\label{cvpier11} \\
		&\delta_k v_{\delta_k} \to 0 \quad \takeshi{\it weakly~in~} \,L^2\takeshi{\bigl(}0,T; H^{3/2}(\Gamma) \takeshi{\bigr)},
	\label{cvpier12}\\
	&\partial_{\boldsymbol{\nu}} u_{\delta_k} -{\delta_k} \Delta_\Gamma v_{\delta_k} \to \partial_{\boldsymbol{\nu}} u \quad \takeshi{\it weakly~in~} 
	\,L^2 (0,T;H_\Gamma).
	\label{cvpier13}
\end{align}
In particular, \takeshi{\eqref{main5}} can be rewritten as
\begin{align}
  \label{pier7}
  w = \partial_\nu u 
  +\eta + \pi_{\Gamma}(v) - g
  \qquad&\text{a.e.~on } \Sigma.
\end{align}
\end{theorem}

\begin{remark}\rm
  We note that the additional assumption \eqref{ip_extra} is a reinforcement of 
  {\rm (A1)} and \eqref{dom2}, for some constant $M \geq \max\{\varrho_1, c_1\}$. 
  In fact, \eqref{ip_extra} implies that the two graphs $\beta$ and $\beta_\Gamma $ have 
  the same growth properties.
\end{remark}

\begin{theorem}
  \label{tak4}
In the setting of Theorem~\ref{tak3}, let $(u,\mu, \xi, v, w, \eta)$ denote the 
  sextuplet solution of the problem~\eqref{main1}--\eqref{main8} given by Theorem~\ref{tak1}  
  and, for $0<\delta\leq1$, let $(u_\delta,\mu_\delta, \xi_\delta, v_\delta, w_\delta, \eta_\delta)$ be the sextuplet solution of the problem~\eqref{weak1}--\eqref{weak7} given by Proposition~\ref{tak0}.  Then,
  there exists a constant $C>0$, independent of $\delta$, such that 
\begin{equation} \label{pier8} 
  \takeshi{|} u_\delta - u\takeshi{|}_{L^\infty(0,T;V^*)\cap L^2(0,T; V)}
  + \takeshi{|}v_\delta - v\takeshi{|}_{L^\infty(0,T;V_\Gamma^*)\cap L^2(0,T; Z_\Gamma)}  \leq C \delta^{1/2}   
\end{equation}
for every $\delta \in (0,1]$ and, as $\delta\to 0$,
\begin{equation} \label{pier9} 
  v_\delta \to v \quad\text{ weakly in } L^2(0,T; V_\Gamma).
\end{equation}
\end{theorem}

\section{Uniform estimates}
\label{uniform}

In this section, we will obtain uniform estimates independent of the 
parameter $0<\delta\leq 1$.
To do so, we recall another suitable approximation to the auxiliary problem. 
Then, taking care of the previous known results, 
we will obtain uniform estimates that are useful for the limiting procedure. 

\subsection{Yosida approximation and viscous Cahn--Hilliard system}

In the approach of~\cite{CFW20}, Proposition~2.1 has been proved by  
considering the following viscous {C}ahn--{H}illiard system: for $\delta, \lambda \in (0,1]$
\begin{align}
  \partial _t u_{\delta, \lambda} - \Delta \mu_{\delta, \lambda} =0 
  & \quad {\rm a.e.\ in ~ } Q, \label{vlw1a}\\
  \mu_{\delta, \lambda} = \lambda \partial_t u_{\delta, \lambda}
  -\Delta u_{\delta, \lambda} + \beta_\lambda(u_{\delta, \lambda}) 
  +\pi(u_{\delta, \lambda})- f
  & \quad {\rm a.e.\ in~} Q, \label{vlw2a}\\
  \partial _{\boldsymbol{\nu}} \mu_{\delta, \lambda} =0 
  & \quad {\rm a.e.\ on~} \Sigma, \label{vlw3a}\\ 
  (u_{\delta, \lambda})_{|_\Gamma} = v_{\delta, \lambda}
  & \quad {\rm a.e.\ on~} \Sigma, \label{vlw4a}\\ 
  \partial _t v_{\delta, \lambda}  - \Delta_\Gamma w_{\delta, \lambda} = 0
  & \quad {\rm a.e.\ on~} \Sigma, \label{vlw5a}\\ 
  w_{\delta, \lambda}= \lambda \partial_t v_{\delta, \lambda}+ \partial_{\boldsymbol{\nu}} u_{\delta, \lambda}- 
  \delta \Delta_\Gamma v_{\delta, \lambda}+\beta_{\Gamma,\lambda } (v_{\delta, \lambda})+\pi_\Gamma(v_{\delta, \lambda})-g
  & \quad {\rm a.e.\ on~} \Sigma, \label{vlw6a}\\
   u_{\delta, \lambda}(0) = u_0 & \quad {\rm a.e.\ in~} \Omega,
  \label{vlw7a} \\
  v_{\delta, \lambda}(0) = v_0 & \quad {\rm a.e.\ on~} \Gamma,
  \label{vlw8a}
\end{align}
where 
$\beta_\lambda $ and $\beta_{\Gamma, \lambda}$ are the {Y}osida approximations of
$\beta $ and $\beta_{\Gamma }$, \luca{respectively,} defined by 
\begin{gather*}
	\beta_\lambda (r):=\frac{1}{\lambda} \bigl( r-J_\lambda (r) \bigr) :=\frac{1}{\lambda} 
	\bigl( r-(I+\lambda \beta) ^{-1}(r)\bigr), \\
	\beta_{\Gamma,\lambda} (r):=\frac{1}{\lambda} \bigl( r-J_{\Gamma,\lambda} (r) \bigr) :=\frac{1}{\lambda} 
	\bigl( r-(I+\lambda \beta_\Gamma) ^{-1}(r)\bigr) \quad {\rm for~} r \in \mathbb{R}.
\end{gather*}
From the well-known theory of maximal monotone operators (see, e.g., \cite{Bar10}), 
we see that $\beta_\lambda$ and $\beta_{\Gamma,\lambda}$ are {L}ipschitz continuous 
functions with {L}ipschitz constant  $1/\lambda$. Moreover, it holds that
\begin{gather}
	\bigl| \beta_\lambda (r) \bigr| \le \bigl|\beta^{\circ}(r) \bigr| ,
	\quad 
	0 \le \widehat{\beta}_\lambda(r) \le \widehat{\beta}(r), 
		\quad {\rm for~all~} r \in D(\beta),
	\label{ap}\\
	\bigl| \beta_{\Gamma,\lambda} (r) \bigr| \le \bigl|\beta_\Gamma^{\circ}(r) \bigr|, 
	\quad 
	0 \le \widehat{\beta}_{\Gamma,\lambda} (r) \le \widehat{\beta}_\Gamma(r)
	\quad {\rm for~all~} r \in D(\beta_\Gamma). 
	\label{apest}
\end{gather}
The approximating problem \eqref{vlw1a}--\eqref{vlw8a} is well posed~\cite{CFW20}, namely, there exists a unique quadruplet  $(u_{\delta, \lambda}, \mu_{\delta, \lambda}, v_{\delta, \lambda}, w_{\delta, \lambda})$, with
\begin{gather*}
	u_{\delta, \lambda} \in H^1(0,T;H) \cap L^\infty (0,T;V) \cap L^2(0,T;W), \\
	\mu_{\delta, \lambda} \in L^2(0,T;W), \\
	v_{\delta, \lambda} \in H^1(0,T;H_\Gamma) \cap L^\infty (0,T;V_\Gamma) \cap L^2(0,T;W_\Gamma), \\
	w_{\delta, \lambda} \in L^2(0,T;W_\Gamma), 
\end{gather*}
satisfying \eqref{vlw1a}--\eqref{vlw8a}. Moreover,
$(u_{\delta,\lambda}, \mu_{\delta,\lambda}, v_{\delta,\lambda}, w_{\delta,\lambda})$ 
converges to the sextuplet $$(u_\delta,\mu_\delta, \xi_\delta, v_\delta, w_\delta, \eta_\delta)$$ given by Proposition~2.1 in a suitable sense, where $\xi_\delta$ and $\eta_\delta$ are the limits of 
$\beta_\lambda(u_{\delta,\lambda})$ and $\beta_{\Gamma,\lambda}(v_{\delta,\lambda})$ as $\lambda \to 0$,
respectively (see, \cite[Theorem~2.3]{CFW20}). 
Therefore, we omit the details of the limiting procedure $\lambda \to 0$ in this paper.

From the next subsection, 
we will obtain the uniform estimates for the approximating problem~\eqref{vlw1a}--\eqref{vlw8a}, whereas
we will discuss the limiting procedure $\delta \to 0$ in the next section. 

\subsection{1st estimate (related to the volume conservation).} 
Integrating \eqref{vlw1a} over $\Omega \times (0,t)  $, 
multiplying by $1/|\Omega|$, and using \eqref{vlw3a}, \eqref{vlw7a} we obtain 
\begin{equation}
	m \bigl(u_{\delta,\lambda}(t) \bigr) =m(u_0) =m_0 
	\label{vc}
\end{equation}
for all $t \in [0,T]$. On the other hand, integrating \eqref{vlw5a} over 
$\Gamma \times (0,t)$ and multiplying by $1/|\Gamma|$, 
from \eqref{vlw8a} we have that 
\begin{equation*}
	m _\Gamma \bigl(v_{\delta,\lambda}(t) \bigr) =m_\Gamma (v_0) =m_{\Gamma0}  
\end{equation*}
for all $t \in [0,T]$. Also, we observe that
\begin{align*}
	\bigl\langle \partial_t \bigl(u_{\delta,\lambda}(t) - m_0 \bigr), 1 \bigr\rangle_{V^*,V} & = \frac{d}{dt} \int_\Omega u_{\delta,\lambda}(t) \, dx 
	 = 0,
\end{align*}
which yields that $\partial_t (u_{\delta,\lambda}(t) - m_0) \in V_0^*$, and analogously 
$\partial_t (v_{\delta,\lambda}(t) - \takeshi{m_{\Gamma0}}) \in V_{\Gamma,0}^*$ for a.a.\ $t \in (0,T)$. 
Moreover, there exists a positive constant $M_1>0$ such that 
\begin{equation}
	\bigl| m(u_{\delta, \lambda } ) \bigr|_{L^\infty(0,T)} 
	+
	\bigl| m_\Gamma(v_{\delta, \lambda } ) \bigr|_{L^\infty(0,T)} 
	\le M_1.
	\label{1st}
\end{equation}

\subsection{2nd estimate}

Multiply \eqref{vlw1a} by $F^{-1}(u_{\delta,\lambda}(t)-u_0)$
and \eqref{vlw5a} by $F^{-1}_\Gamma (v_{\delta,\lambda}(t)-v_0)$. 
Then,
using \eqref{vlw3a} we obtain 
\begin{equation}
	\bigl \langle \partial_t \bigl(u_{\delta,\lambda}(t)-u_0 \bigr), F^{-1} \bigl(u_{\delta,\lambda}(t)-u_0 \bigr) \bigr\rangle_{V_0^*,V_0}
	+ \int_\Omega \nabla \mu_{\delta,\lambda}(t) \cdot \nabla F^{-1} \bigl( u_{\delta,\lambda}(t)-u_0 \bigr)
	\,dx =0,
	\label{2nd1}
\end{equation}
and 
\begin{equation}
	\bigl \langle \partial_t \bigl(v_{\delta,\lambda}(t)-v_0 \bigr), 
	F^{-1}_\Gamma \bigl(v_{\delta,\lambda}(t)-v_0 \bigr) \bigr\rangle_{V_{\Gamma,0}^*,V_{\Gamma,0} }
	\! + \! \int_\Gamma \nabla_\Gamma w_{\delta,\lambda}(t) \cdot 
	\nabla_\Gamma F^{-1}_\Gamma \bigl( v_{\delta,\lambda}(t)-v_0 \bigr)\,d\Gamma =0
	\label{2nd4}
\end{equation}
for a.a.\ $t \in (0,T)$. 
Next, multiplying \eqref{vlw2a} by $u_{\delta,\lambda}(t)-u_0$
and using \eqref{vlw4a} we infer that
\begin{align}
	& \bigl( \mu_{\delta,\lambda} (t), u_{\delta,\lambda}(t)-u_0\bigr)_{\!H} 
	\nonumber \\
	& = \frac{\lambda}{2} \frac{d}{dt} \bigl| u_{\delta,\lambda}(t)-u_0 \bigr|_{H_0}^2 + 
	\int _\Omega \nabla u_{\delta, \lambda}(t)\cdot\nabla(u_{\delta, \lambda}(t)-u_0)\,dx
	- \bigl(\partial_{\boldsymbol{\nu}} u_{\delta, \lambda}(t) ,  
	v_{\delta, \lambda}(t)  -v_0 \bigr)_{\!H_\Gamma}
	\nonumber \\
	& \quad {}
	+ \bigl( \beta_{\lambda} \bigl( u_{\delta, \lambda}(t)  \bigr), u_{\delta, \lambda}(t) -u_0\bigr)_{\!H}
	+ \bigl( \pi \bigl( u_{\delta, \lambda}(t)  \bigr)-f(t), u_{\delta, \lambda}(t) -u_0\bigr)_{\!H}.
	\label{2nd2}
\end{align}
Analogously, multiplying \eqref{vlw6a} by $v_{\delta,\lambda}(t)-v_0$ we have that
\begin{align}
	& \bigl( w_{\delta,\lambda} (t), v_{\delta,\lambda}(t)-v_0\bigr)_{\!H_\Gamma} 
	\nonumber \\
	& = 
	\frac{\lambda}{2} \frac{d}{dt} \bigl| v_{\delta,\lambda}(t)-v_0 \bigr|_{H_{\Gamma,0}}^2 
	+
	\bigl(\partial_{\boldsymbol{\nu}} u_{\delta, \lambda}(t) , 
	v_{\delta, \lambda}(t) -v_0\bigr)_{\!H_\Gamma}
	+ \delta
	\int _\Gamma \nabla_\Gamma v_{\delta, \lambda}(t) 
	\cdot\nabla_\Gamma(v_{\delta, \lambda}(t)-v_0)\,d\Gamma 
	\nonumber \\
	& \quad  {} + \bigl( \beta_{\Gamma,\lambda} 
	\bigl( v_{\delta, \lambda}(t)  \bigr), v_{\delta, \lambda}(t) -v_0\bigr)_{\!H_\Gamma}
	+ \bigl( \pi_\Gamma \bigl( v_{\delta, \lambda}(t)  \bigr)-g(t), 
	v_{\delta, \lambda}(t) -v_0\bigr)_{\!H_\Gamma}.
	\label{2nd3}
\end{align}
By merging \eqref{2nd1}, \eqref{2nd4}, \eqref{2nd2}, and \eqref{2nd3}, and then
adding $|u_{\delta,\lambda}(t)|_H^2$ to both sides of the resultant we obtain 
\begin{align}
	& \frac{1}{2} \frac{d}{dt} \bigl| u_{\delta,\lambda}(t)-u_0 \bigr|_{V_0^*}^2
	+ 
	\frac{\lambda}{2} \frac{d}{dt} \bigl| u_{\delta,\lambda}(t)-v_0 \bigr|_{H_0}^2 
	+ 
	\frac{1}{2} \frac{d}{dt} \bigl| v_{\delta,\lambda}(t)-v_0\bigr|_{V_{\Gamma,0}^*}^2 
	\nonumber \\
	& \quad {} 
	+
	\frac{\lambda}{2} \frac{d}{dt} \bigl| v_{\delta,\lambda}(t)-v_0 \bigr|_{H_{\Gamma,0}}^2 
	+
	\bigl| u_{\delta,\lambda}(t) \bigr|_{V}^2
	+
	\delta
	\int _\Gamma \bigl| \nabla_\Gamma v_{\delta, \lambda}(t) \bigr|^2\, d\Gamma 
	\nonumber \\
	& \quad {} +
	\bigl( \beta_{\lambda} \bigl( u_{\delta, \lambda}(t)  \bigr), u_{\delta, \lambda}(t) -u_0 \bigr)_{\!H}
	+ \bigl( \beta_{\Gamma,\lambda} \bigl( v_{\delta, \lambda}(t)  \bigr), 
	v_{\delta, \lambda}(t) -v_0 \bigr)_{\!H_\Gamma}
	\nonumber \\
	& \le \bigl| u_{\delta,\lambda}(t) \bigr|_H^2 +
	\int_\Omega\nabla u_{\delta,\lambda}(t)\cdot\nabla u_0\,dx
	+\delta \int_\Gamma\nabla_\Gamma v_{\delta,\lambda}(t)\cdot\nabla_\Gamma v_0\,d\Gamma
	\nonumber \\
	& \quad {} 
	- \bigl( \pi \bigl( u_{\delta, \lambda}(t)  \bigr)-f(t), u_{\delta, \lambda}(t) -u_0\bigr)_{\!H}
	- \bigl( \pi_\Gamma \bigl( v_{\delta, \lambda}(t)  \bigr)-g(t), v_{\delta, \lambda}(t) -v_0\bigr)_{\!H_\Gamma}
	\label{2nd5}
\end{align}
for a.a.\ $t \in (0,T)$. Now, on the left-hand side,
by the convexity of $\widehat\beta_\lambda$ and $\widehat\beta_{\Gamma,\lambda}$,
as well as \eqref{ap}--\eqref{apest}, we deduce that
\begin{align*}
 &\bigl( \beta_{\lambda} \bigl( u_{\delta, \lambda}(t)  \bigr), u_{\delta, \lambda}(t) -u_0 \bigr)_{\!H}
	+ \bigl( \beta_{\Gamma,\lambda} \bigl( v_{\delta, \lambda}(t)  \bigr), 
	v_{\delta, \lambda}(t) -v_0 \bigr)_{\!H_\Gamma}\\
 &\geq
  \int_\Omega \widehat{\beta}_{\lambda} \bigl( u_{\delta, \lambda}(t)  \bigr) \,dx
  - \int_\Omega \widehat{\beta} ( u_0 ) \,dx  
 +  \int_\Gamma \widehat{\beta}_{\Gamma, \lambda} \bigl( v_{\delta, \lambda}(t) \bigr) \,d\Gamma 
 -\int_\Gamma \widehat{\beta}_{\Gamma} ( v_0 ) \,d\Gamma.
\end{align*}
On the right-hand side, by the Young inequality we have
\begin{align*}
  &\int_\Omega\nabla u_{\delta,\lambda}(t)\cdot\nabla u_0\,dx
  +\delta \int_\Gamma\nabla_\Gamma v_{\delta,\lambda}(t)\cdot\nabla_\Gamma v_0\,d\Gamma\\
  &\leq \frac12|u_{\delta,\lambda}(t)|_V^2 
  + \frac\delta2\int _\Gamma \bigl| \nabla_\Gamma v_{\delta, \lambda}(t) \bigr|^2 \, d\Gamma
  + |u_0|_V^2 + \delta|v_0|^2_{V_\Gamma}.
\end{align*}
Furthermore, 
applying the {E}hrling lemma (see, e.g., \cite[Chapter~1, Lemma~5.1]{Lio69}) for  
$V \mathop{\hookrightarrow} \mathop{\hookrightarrow} H \subset V^*$, 
we see that for any $\varepsilon>0$ there 
exists a positive constant $C_\varepsilon>0$ such that
\begin{equation}
	\bigl| u_{\delta,\lambda}(t) \bigr|_{H}^2
	 \le \varepsilon \bigl| u_{\delta,\lambda}(t)\bigr|_{V}^2 
	+ C_\varepsilon  \Bigl(1+\bigl| u_{\delta,\lambda}(t) -u_0\bigr|_{V_0^*}^2 \Bigr),
	\label{pier1}
\end{equation}
where we have added and subtracted $u_0$ in the second term on the right-hand side
and used the equivalence of $|\cdot|_{V^*}$ and $|\cdot|_{V_0^*}$ on $V_0^*$.
Moreover, thanks to {\rm (A2)} and {\rm (A3)},
using the Young inequality and again the {E}hrling lemma 
it turns out that
there exists a positive constant $\takeshi{C>0}$ depending on 
$\pi$, $|u_0|_H$, and $|\Omega|$ such that 
\begin{align*}
	& -\bigl( \pi \bigl( u_{\delta, \lambda}(t)  \bigr)-f(t), u_{\delta, \lambda}(t) -u_0\bigr)_{\!H}
	\\
	& \le \bigl( L \bigl| u_{\delta, \lambda}(t)  \bigr|_H + |\pi(0)|_H + |f(t)|_H 
	\bigr) \bigl(\bigl| u_{\delta, \lambda}(t)  \bigr|_H	
	 + |u_0|_H \bigr)
	\\
	& \le \varepsilon \bigl| u_{\delta,\lambda}(t) \bigr|_{V}^2 
	+ C_\varepsilon \Bigl(1+ \bigl| u_{\delta,\lambda}(t) - u_0\bigr|_{V_0^*}^2\Bigr)
	+ \takeshi{C} \left(1+|f(t)|_H^2\right)
\end{align*}
for a.a.\ $t \in (0,T)$. Analogously,  
using the Young inequality, the {E}hrling lemma 
with respect to the inclusions 
$Z_{\Gamma} \mathop{\hookrightarrow} \mathop{\hookrightarrow}  H_\Gamma \mathop{\hookrightarrow}V_{\Gamma}^*$, and the estimate \eqref{tr0} for the trace $\gamma_0$,
we deduce that
\begin{align*}
	&-\bigl( \pi_\Gamma \bigl( v_{\delta, \lambda}(t)  \bigr)-g(t), 
	v_{\delta, \lambda}(t) -v_0\bigr)_{\!H_\Gamma}\\
	& \le \bigl( L_\Gamma \bigl| v_{\delta, \lambda}(t)  \bigr|_{H_\Gamma} 
	+ |\pi_\Gamma(0)|_{H_\Gamma} + |g(t)|_{H_\Gamma}
	\bigr) \bigl(\bigl| v_{\delta, \lambda}(t)  \bigr|_{H_\Gamma}
	 + |v_0|_{H_\Gamma} \bigr)
	 \\
	&\le 	
	\varepsilon \bigl| u_{\delta,\lambda}(t) \bigr|_{V}^2 
	+ C_\varepsilon \bigl| v_{\delta,\lambda}(t) -v_0\bigr|_{V_{\Gamma,0}^*}^2
	+ \takeshi{C} \left(1+|g(t)|_{H_\Gamma}^2\right)
\end{align*}
for a.a.\ $t \in (0,T)$, where
we exploited the equivalence of $|\cdot|_{V_\Gamma^*}$ and $|\cdot|_{V_{\Gamma,0}^*}$ 
on $V_{\Gamma,0}^*$ and
we let the updated value of \takeshi{$C$}
depend also on 
$\pi_\Gamma$, $|v_0|_{H_\Gamma}$, and $|\Gamma|$.
Therefore, going back to \eqref{2nd5}
we choose $\varepsilon$ small enough and
apply the {G}ronwall inequality, obtaining
\begin{align}
	\nonumber
	& \sup_{t\in[0,T]}\bigl| u_{\delta,\lambda}(t)-u_0 \bigr|_{V_0^*}^2
	+ 
	\sup_{t\in[0,T]}\lambda\bigl| u_{\delta,\lambda}(t)-u_0 \bigr|_{H_0}^2 \\
	& \quad {} + 
	\sup_{t\in[0,T]}\bigl| v_{\delta,\lambda}(t)-v_0 \bigr|_{V_{\Gamma,0}^*}^2 
	+
	\sup_{t\in[0,T]}\lambda \bigl| v_{\delta,\lambda}(t)-v_0 \bigr|_{H_{\Gamma,0}}^2 
	\nonumber \\
	& \quad {} 
	+
	\int_0^{T} \bigl| u_{\delta,\lambda}(s) \bigr|_{V}^2\,ds 
	+
	\delta
	\int_0^{T} 
	\bigl| \nabla_\Gamma v_{\delta, \lambda}(s) \bigr|_{H_\Gamma}^2 \,ds 
	\nonumber
	\\
	& \quad {}
	+ \int_0^T \bigl| \widehat{\beta}_\lambda \bigl(u_{\delta, \lambda}(s)  \bigr) \bigr|_{L^1(\Omega)} \,ds
	+ \int_0^T \bigl| \widehat{\beta}_{\Gamma,\lambda} 
	\bigl(v_{\delta, \lambda}(s)  \bigr) \bigr|_{L^1(\Gamma)}\, ds 
	\le M_2,
	\label{2ndf}
\end{align}
where the constant $M_2$ depends on \takeshi{$T$}, $|f|_{L^2(0,T; H)}$,
$|g|_{L^2(0,T; H_\Gamma)}$, $|u_0|_V$, and $\delta^{1/2}|v_0|_{V_\Gamma}$.

\subsection{3rd estimates}

Firstly, multiplying \eqref{vlw1a} by $\mu_{\delta,\lambda}(t)+f(t)$ and using \eqref{vlw3a} we obtain 
\begin{align}
	\nonumber
	&\bigl \langle \partial_t u_{\delta,\lambda}(t), \mu_{\delta,\lambda}(t)+f(t) \bigr\rangle_{V^*,V}
	+ \int_\Omega \bigl| \nabla \mu_{\delta,\lambda}(t) \bigr|^2\, dx \\
	& = -\int_\Omega \nabla \mu_{\delta,\lambda}(t) \cdot \nabla f(t)\, dx
	 \le 
	\frac{1}{2}\int_\Omega \bigl| \nabla \mu_{\delta,\lambda}(t) \bigr|^2\, dx +
	\frac{1}{2}\int_\Omega \bigl| \nabla f(t) \bigr|^2\, dx
	\label{3rd1}
\end{align}
for a.a.\ $t \in (0,T)$. 
Secondly, multiplying \eqref{vlw2a} by $\partial_t u_{\delta,\lambda}(t)$ leads to
\begin{align}
	& \bigl \langle \partial_t u_{\delta,\lambda}(t), \mu_{\delta,\lambda}(t)+f(t) \bigr\rangle_{V^*,V}
	\nonumber \\
	& = \lambda \bigl| \partial_t u_{\delta,\lambda}(t) \bigr|_{H}^2 + 
	\frac{1}{2}	
	\frac{d}{dt}
	\int _\Omega \bigl| \nabla u_{\delta, \lambda}(t) \bigr|^2\,dx 
	- \bigl(\partial_{\boldsymbol{\nu}} u_{\delta, \lambda}(t),  
	 \partial _t v_{\delta, \lambda}(t)\bigr)_{\!H_\Gamma}
	\nonumber \\
	& \quad {}
	+ \frac{d}{dt} \left\{ \int_\Omega \widehat{\beta}_{\lambda} \bigl( u_{\delta, \lambda}(t)\bigr)\, dx
	+
	\int_\Omega \widehat{\pi} \bigl( u_{\delta, \lambda}(t)  \bigr) \,dx \right\}.
	\label{3rd2}
\end{align}
Analogously, multiplying \eqref{vlw6a} by $\partial_t v_{\delta,\lambda}(t)$ we infer that
\begin{align}
	& \bigl \langle \partial_t v_{\delta,\lambda}(t), 
	w_{\delta,\lambda}(t) +g(t) \bigr\rangle_{V_\Gamma^*,V_\Gamma}
	\nonumber \\
	& = \lambda \bigl| \partial_t v_{\delta,\lambda}(t) \bigr|_{H_\Gamma}^2 
	+ \bigl(\partial_{\boldsymbol{\nu}} u_{\delta, \lambda}(t) ,  
	\partial _t v_{\delta, \lambda}(t) \bigr)_{\!H_\Gamma}
	+
	\frac{\delta}{2}
	\frac{d}{dt}
	\int _\Gamma \bigl| \nabla_\Gamma v_{\delta, \lambda}(t) \bigr|^2\,d\Gamma 
	\nonumber \\
	& \quad {}
	+ \frac{d}{dt} \left\{ \int_\Gamma \widehat{\beta}_{\Gamma,\lambda} 
	\bigl( v_{\delta, \lambda}(t) \bigr)\, d\Gamma
	+
	\int_\Gamma \widehat{\pi}_\Gamma \bigl( v_{\delta, \lambda}(t)  \bigr) \,d\Gamma \right\},
	\label{3rd3}
\end{align}
while multiplying \eqref{vlw5a}  by $w_{\delta,\lambda}(t)+g(t)$ gives
\begin{align}
	\nonumber
	&\bigl \langle \partial_t v_{\delta,\lambda}(t), 
	w_{\delta,\lambda}(t)+g(t) \bigr\rangle_{V_{\Gamma}^*,V_{\Gamma} }
	+ \int_\Gamma \bigl| \nabla_\Gamma w_{\delta,\lambda}(t) \bigr|^2 \,d\Gamma \\
	&\le \frac{1}{2} \int_\Gamma \bigl| \nabla_\Gamma w_{\delta,\lambda}(t) \bigr|^2 \, d\Gamma
	+
	\frac{1}{2} \int_\Gamma \bigl| \nabla_\Gamma g(t) \bigr|^2 \, d\Gamma.
	\label{3rd4}
\end{align}
Combining \eqref{3rd1}--\eqref{3rd4},
integrating the resulting inequality from $0$ to $t$, 
adding the term $(1/2)|u_{\delta,\lambda}(t)|_H^2$ to both sides,
and using \eqref{vlw7a}--\eqref{apest}, we obtain 
\begin{align}
	& \frac{1}{2} \int_0^t \bigl| \nabla \mu_{\delta,\lambda}(s) \bigr|_{H}^2 \,ds
	+ 
	\frac{1}{2} \int_0^t \bigl| \nabla_\Gamma w_{\delta,\lambda}(s) \bigr|_{H_\Gamma}^2\, ds
	\nonumber \\
	& \quad {}+
	\lambda \int_0^t \bigl| \partial _t u_{\delta,\lambda}(s) \bigr|_{H}^2 \,ds
	+ 
	\lambda \int_0^t \bigl| \partial _t v_{\delta,\lambda}(s) \bigr|_{H_\Gamma}^2 \,ds
	+
	\frac{1}{2} \bigl| u_{\delta,\lambda}(t) \bigr|_{V}^2 \nonumber \\
	& \quad {} 
	+
	\frac{\delta}{2}
	\bigl| \nabla_\Gamma v_{\delta,\lambda}(t) \bigr|_{H_\Gamma}^2
	+\int_\Omega \widehat{\beta}_\lambda \bigl( u_{\delta,\lambda}(t) \bigr) \,dx 
	+\int_\Gamma \widehat{\beta}_{\Gamma,\lambda} \bigl( v_{\delta,\lambda}(t) \bigr) \,d\Gamma
	\nonumber \\
	\nonumber
	& \le 
	\frac{1}{2} | \nabla u_0 |_{H}^2 
	+
	\frac{\delta}{2}
	| \nabla_\Gamma v_0 |_{H_\Gamma}^2
	+\int_\Omega \widehat{\beta} ( u_0 ) \,dx 
	+\int_\Gamma \widehat{\beta}_{\Gamma} ( v_0 ) \,d\Gamma \\
	\nonumber
	&\quad {}
	+\frac{1}{2} \int_0^T \bigl| f(s) \bigr|_V^2\, ds 
	+ \frac{1}{2} \int_0^T \bigl| g(s) \bigr|_{V_\Gamma}^2 \,ds
	+ \frac{1}{2} \bigl|u_{\delta,\lambda}(t) \bigr|_H^2\\
	&\quad{}
	+\int_\Omega \bigl| \widehat{\pi} \bigl( u_{\delta,\lambda}(t) \bigr) \bigr| \,dx 
	+\int_\Omega \bigl| \widehat{\pi}(u_0) \bigr| \,dx 
	+\int_\Gamma \bigl| \widehat{\pi}_\Gamma \bigl( v_{\delta,\lambda}(t) \bigr) \bigr| \,d\Gamma
	+\int_\Gamma \bigl| \widehat{\pi}_\Gamma(v_0) \bigr| \,d\Gamma 
	\label{3rd6}
\end{align}
for all $t \in [0,T]$. 
Here, from {\rm (A2)} we see that 
\begin{align*}
	\bigl| \widehat{\pi}(r) \bigr| \le L |r|^2+ \frac{1}{2L} \bigl| \pi(0) \bigr|^2,
	\qquad
	\bigl| \widehat{\pi}_\Gamma(r) \bigr| \le L_\Gamma |r|^2+ \frac{1}{2L_\Gamma} \bigl| \pi_\Gamma(0) \bigr|^2
\end{align*}
for all $r \in \mathbb{R}$, therefore
\begin{align*}
	\int_\Omega \bigl| \widehat{\pi} \bigl( u_{\delta,\lambda}(t) \bigr) \bigr| \,dx 
	+
	\int_\Omega \bigl| \widehat{\pi}(u_0) \bigr| \, dx 
	&\le L\bigl|u_{\delta,\lambda}(t) \bigr|_H^2+ L |u_0|_H^2+\frac{1}{L} \bigl| \pi(0) \bigr|^2, 
	\nonumber
	\\
	\int_\Gamma \bigl| \widehat{\pi}_\Gamma \bigl( v_{\delta,\lambda}(t) \bigr) \bigr| \, d\Gamma
	+
	\int_\Gamma \bigl| \widehat{\pi}_\Gamma(v_0) \bigr| \, d\Gamma 
	&\le L_\Gamma \bigl|v_{\delta,\lambda}(t) \bigr|_{H_\Gamma}^2
	+L_\Gamma |v_0|_{H_\Gamma}^2+\frac{1}{L_\Gamma} \bigl| \pi_\Gamma(0) \bigr|^2.
\end{align*}
Now, 
applying again the compactness inequalities and the estimate~\eqref{tr0} for the 
trace, we see that for any $\varepsilon>0$
there exists a constant $C_\varepsilon>0$ such that \eqref{pier1} and 
\begin{equation}
	\bigl|v_{\delta,\lambda}(t) \bigr|_{H_\Gamma}^2 
	\le \varepsilon \bigl|u_{\delta,\lambda}(t) \bigr|_V^2 
	+ C_\varepsilon \Bigl( 1+ \bigl|v_{\delta,\lambda}(t) -v_0 \bigr|_{V_{\Gamma,0}^*}^2 
	\Bigr)
	\label{pier2}
\end{equation}
hold, where $C_\varepsilon$ depends on $|u_0|_H$, $|v_0|_{H_\Gamma}$, $|\Omega|$, and 
$|\Gamma|$. Thus, using \eqref{2ndf}, we deduce that there exists a positive constant $M_3>0$ such that 
\begin{align}
	& \int_0^T \bigl| \nabla \mu_{\delta,\lambda}(s) \bigr|_{H}^2 \,ds
	+ 
	\int_0^T \bigl| \nabla_\Gamma w_{\delta,\lambda}(s) \bigr|_{H_\Gamma}^2 \,ds
	\nonumber \\
	& \quad {} 
	+
	\lambda \int_0^T \bigl| \partial _t u_{\delta,\lambda}(s) \bigr|_{H}^2 \,ds
	+ 
	\lambda \int_0^T \bigl| \partial _t v_{\delta,\lambda}(s) \bigr|_{H_\Gamma}^2 \,ds
	\nonumber \\
	& \quad {} 
	+\sup_{t\in[0,T]}
	\bigl| u_{\delta,\lambda}(t) \bigr|_{V}^2 
	+\sup_{t\in[0,T]}
	\delta
	\bigl| \nabla_\Gamma v_{\delta,\lambda}(t) \bigr|_{H_\Gamma}^2
	\nonumber \\
	& \quad {} 
	+\sup_{t\in[0,T]}\int_\Omega \widehat{\beta}_\lambda \bigl( u_{\delta,\lambda}(t) \bigr) \,dx 
	+\sup_{t\in[0,T]}\int_\Gamma \widehat{\beta}_{\Gamma,\lambda} 
	\bigl( v_{\delta,\lambda}(t) \bigr) \,d\Gamma \le M_3.
	\label{3rdf}
\end{align}
From \eqref{vlw1a}, \eqref{vlw3a},  and \eqref{vlw5a}, 
it straightforward to infer that
\begin{equation*}
	\bigl| \partial _t u_{\delta,\lambda}(s) \bigr|_{V_0^*}^2 
	\le \bigl| \nabla \mu_{\delta,\lambda}(s) \bigr|_{H}^2, \quad 
	\bigl| \partial _t v_{\delta,\lambda}(s) \bigr|_{V_{\Gamma,0}^*}^2 
	\le \bigl| \nabla_\Gamma w_{\delta,\lambda}(s) \bigr|_{H_\Gamma}^2, 
\end{equation*}
for a.a.\ $s \in (0,T)$.
Thus, the estimate \eqref{3rdf} also implies that
\begin{equation}
	\int_0^T \bigl| \partial _t u_{\delta,\lambda}(s) \bigr|_{V_0^*}^2 \,ds
	+ 
	\int_0^T \bigl| \partial _t v_{\delta,\lambda}(s) \bigr|_{V_{\Gamma,0}^*}^2 \,ds
	\le M_3 .
	\label{3rdf2}
\end{equation}

\subsection{4th estimate.}

Thanks to {\rm (A1)} and {\rm (A4)}, we can use the following 
useful inequality (see \cite[Appendix, Prop.\ A.1]{MiZe} and/or 
the detailed proof given in \cite[p.\ 908]{GMS09}): 
there exist two positive constants $c_2, c_3>0$ such that
\begin{gather*}
	\bigl( \beta_\lambda \bigl( u_{\delta,\lambda}(t) \bigr), u_{\delta,\lambda}(t) -u_0 \bigr)_{\!H} 
	\ge c_2\int_\Omega \bigl| \beta_\lambda \bigl( u_{\delta,\lambda}(t) \bigr) \bigr|\,dx -c_3 |\Omega|, \\
	\bigl( \beta_{\Gamma,\lambda} \bigl( v_{\delta,\lambda}(t) \bigr), 
	v_{\delta,\lambda}(t) -v_{0} \bigr)_{\!H_\Gamma} 
	\ge c_2 \int_\Gamma \bigl| \beta_{\Gamma,\lambda} 
	\bigl( v_{\delta,\lambda}(t) \bigr) \bigr|\,d\Gamma -c_3 |\Gamma|
\end{gather*}
for a.a.\ $t \in(0,T)$. Therefore, merging \eqref{2nd1}--\eqref{2nd3} again
and recalling the definition of inner products of $V_0^*$ and $V_{\Gamma,0}^*$, we have 
\begin{align}
	& c_2 \left\{ \int_\Omega \bigl| \beta_\lambda \bigl( u_{\delta,\lambda}(t) \bigr) \bigr| \,dx 
	+
	\int_\Gamma \bigl| \beta_{\Gamma,\lambda} \bigl( v_{\delta,\lambda}(t) \bigr) \bigr| \,d\Gamma 
	\right\} 
	- c_3 \bigl( |\Omega| + |\Gamma| \bigr) 
	\nonumber \\
	& \le
	\bigl( f(t)-\pi\bigl( u_{\delta,\lambda}(t) \bigr)
	-\lambda \partial_t  u_{\delta,\lambda}(t), u_{\delta,\lambda}(t) -u_0 \bigr)_{\!H} 
	- \bigl( \partial_t  u_{\delta,\lambda}(t), u_{\delta,\lambda}(t) -u_0 \bigr)_{V_0^*}
	\nonumber \\
	& \quad {}
	+ 
	\bigl( g(t)-\pi_\Gamma \bigl( v_{\delta,\lambda}(t) \bigr)
	-\lambda \partial_t  v_{\delta,\lambda}(t), v_{\delta,\lambda}(t) -v_0 \bigr)_{\!H_\Gamma} 
	- \bigl( \partial_t  v_{\delta,\lambda}(t), v_{\delta,\lambda}(t) -v_0 \bigr)_{V_{\Gamma,0}^*}
	 \nonumber \\
	 & \le 
	 \Bigl\{ 
	 \bigl| f(t) \bigr|_H  + \bigl| \pi\bigl( u_{\delta,\lambda}(t) \bigr) \bigr|_H
	 + \lambda \bigl| \partial_t  u_{\delta,\lambda}(t) \bigr|_H
	\Bigl\} 	 
	 \bigl| u_{\delta,\lambda}(t) -u_0 \bigr|_H 
	\nonumber \\
	& \quad {}
	+ 
	\Bigl\{ 
	\bigl| g(t) \bigr|_{H_\Gamma}  
	+ \bigl| \pi_\Gamma \bigl( v_{\delta,\lambda}(t) \bigr) \bigr|_{H_\Gamma}
	 + \lambda \bigl| \partial_t  v_{\delta,\lambda}(t) \bigr|_{H_\Gamma} 
	 \Bigr\} 
	 \bigl| v_{\delta,\lambda}(t) -v_0 \bigr|_{H_\Gamma} 
	 \nonumber \\
	& \quad {}
	+ \bigl| \partial_t  u_{\delta,\lambda}(t) \bigr|_{V_0^*} \bigl| u_{\delta,\lambda}(t) -u_0 \bigr|_{V_0^*}
	 + \bigl| \partial_t  v_{\delta,\lambda}(t) \bigr|_{V_{\Gamma,0}^*} 
	 \bigl| v_{\delta,\lambda}(t) -v_0 \bigr|_{V_{\Gamma,0}^*}
	 \label{4th1}
\end{align}
for a.a.\ $t \in (0,T)$. Here, from {\rm (A2)} and \eqref{3rdf}--\eqref{3rdf2} 
it follows that the right-hand side of \eqref{4th1} is uniformly bounded in $L^2(0,T)$:
hence, 
there exists a positive constant $M_4>0$ such that 
\begin{equation}
	\int_0^T \bigl| \beta_\lambda \bigl( u_{\delta,\lambda}(s) \bigr) \bigr|_{L^1(\Omega)}^2 \,ds 
	+
	\int_0^T \bigl| \beta_{\Gamma,\lambda} \bigl( v_{\delta,\lambda}(s) \bigr) \bigr|_{L^1(\Gamma)}^2 \,ds
	\le M_4.
	\label{4thf}
\end{equation}

\subsection{5th estimate}
Setting $W_0:=H^2(\Omega)\cap H^1_0(\Omega)$, 
we multiply \eqref{vlw2a} by an arbitrary $\takeshi{\zeta} \in L^2(0,T; W_0)$ and
integrate by parts. Recalling the continuous embedding
$W_0\mathop{\hookrightarrow} L^\infty(\Omega)$,
we obtain that 
\begin{align*}
  & \int_0^T \takeshi{\bigl(} \mu_{\delta,\lambda}(s), \takeshi{\zeta}(s) 
  \takeshi{\bigr)_{\!H}}\,ds \\
  & \leq
  \int_0^T\bigl|\lambda\partial_t u_{\delta,\lambda}(s)
  +\pi(u_{\delta,\lambda}(s))-f(s)\bigr|_H \takeshi{\bigl|} 
  \takeshi{\zeta}(s) \takeshi{\bigr|}_H\,ds\\
  &\quad{}+\int_0^T\bigl|\nabla u_{\delta,\lambda}(s)\bigr|_H
  \bigl|\nabla \takeshi{\zeta}(s)\bigr|_H\,ds +C\int_0^T\bigl|\beta_{\lambda} 
  \takeshi{\bigl(} u_{\delta,\lambda}(s) \takeshi{\bigr)} \bigr|_{L^1(\Omega)}
  \takeshi{\bigl|} \takeshi{\zeta}(s) \takeshi{\bigr|}_{W_0}\,ds,
\end{align*}
where the positive constant $C$ only depends on $\Omega$.
Therefore, exploiting the estimates \eqref{3rdf} and \eqref{4thf} we infer that 
\begin{equation}
	\int_0^T \bigl| \mu_{\delta,\lambda}(s)\bigr|_{W_0^*}^2 \,ds 
	\le M_5.
	\label{5pier}
\end{equation}
Now, we apply the Ehrling lemma for the spaces 
$V\mathop{\hookrightarrow}\mathop{\hookrightarrow} H \mathop{\hookrightarrow} W_0^*$
to deduce that for every $\varepsilon>0$ there exists a constant $C_\varepsilon>0$ such that 
\[
  \bigl| \mu_{\delta,\lambda}(s)\bigr|_{H}^2
  \leq \varepsilon \bigl| \nabla\mu_{\delta,\lambda}(s)\bigr|_{H}^2
  +C_\varepsilon \bigl| \mu_{\delta,\lambda}(s)\bigr|_{W_0^*}^2
  \quad\text{for a.a.~$s\in(0,T)$.}
\]
Consequently, the estimates \eqref{3rdf} and \eqref{5pier} yield, possibly updating $M_5$,
\begin{equation}
	\int_0^T \bigl| \mu_{\delta,\lambda}(s)\bigr|_{V}^2 \,ds 
	\le M_5.
	\label{6pier}
\end{equation}
Next, we test~\eqref{vlw2a} by $1$ and integrate by parts using the boundary equations \eqref{vlw4a} and \eqref{vlw6a}. Recalling that $\partial_t u_{\delta,\lambda}(t) \in V_0^*$ and
$\partial_t v_{\delta,\lambda}(t) \in V_{\Gamma,0}^*$, it easily follows that 
\begin{align}
	& \int_\Omega \mu_{\delta,\lambda} (t) \, dx  + \int_\Gamma w_{\delta,\lambda} (t) \, d\Gamma	\nonumber \\
	& =  \int_\Omega \beta_{\lambda} \bigl( u_{\delta, \lambda}(t)  \bigr) \, dx 
	+ \bigl( \pi \bigl( u_{\delta, \lambda}(t)  \bigr)-f(t), 1 \bigr)_{\!H}
	\nonumber \\
	& \quad  {} + \int_\Gamma \beta_{\Gamma,\lambda} 
	\bigl( v_{\delta, \lambda}(t)  \bigr)\, \takeshi{d\Gamma}
	+ \bigl( \pi_\Gamma \bigl( v_{\delta, \lambda}(t)  \bigr)-g(t), 1\bigr)_{\!H_\Gamma}	\label{7pier}
\end{align}
for a.a.\ $t \in (0,T)$. Then, by virtue of \eqref{2ndf}, \eqref{4thf}, \eqref{6pier} 
and assumptions~(A2) and (A3), comparing the terms in \eqref{7pier} yields that 
$$
\hbox{the function } \ t \mapsto
 m_\Gamma (w_{\delta,\lambda} (t) ) = \frac1{|\Gamma|} \int_\Gamma w_{\delta,\lambda} (t) \, d\Gamma \ \hbox{ is uniformly bounded in $L^2(0,T)$,}
$$
whence the estimate \eqref{3rdf} and the Poincar\'e--Wirtinger inequality allow us to infer~that
\begin{equation}
	\int_0^T \bigl| w_{\delta,\lambda}(s)\bigr|_{V_\Gamma}^2 \,ds 
	\le M_5.
	\label{1luca}
\end{equation}

\subsection{6th and 7th estimates}
We test now equation \eqref{vlw2a} by $\beta_\lambda(u_{\delta,\lambda}
)$ and
equation \eqref{vlw6a} by $\beta_\lambda(v_{\delta,\lambda}
)$, then we combine them obtaining
\begin{align}
  &\frac\lambda2\frac{d}{dt}\int_\Omega\widehat\beta_\lambda(u_{\delta,\lambda})\,dx +
  \int_\Omega\beta_\lambda'(u_{\delta,\lambda})|\nabla u_{\delta,\lambda}|^2\,dx +
  \int_\Omega |\beta_\lambda(u_{\delta,\lambda})|^2\,dx \nonumber \\
  &\qquad+
  \frac\lambda2\frac{d}{dt}\int_\Gamma\widehat\beta_\lambda(v_{\delta,\lambda})\,d\Gamma +
  \delta\int_\Gamma\beta_\lambda'(v_{\delta,\lambda})|\nabla_\Gamma v_{\delta,\lambda}|^2\,d\Gamma +
  \int_\Gamma \beta_\lambda(v_{\delta,\lambda})\beta_{\Gamma,\lambda}(v_{\delta,\lambda})\,d\Gamma \nonumber \\
  &=\int_\Omega\left(\mu_{\delta,\lambda} + f - \pi(u_{\delta,\lambda})\right)
  \beta_\lambda(u_{\delta,\lambda})\,dx + 
  \int_\Gamma\left(w_{\delta,\lambda} + g - \pi_\Gamma(v_{\delta,\lambda})\right)
  \beta_\lambda(v_{\delta,\lambda})\,d\Gamma \label{8pier}
\end{align}
Now, we recall assumption~(A1) and point out that \eqref{dom2} entails that the same 
inequality holds for the Yosida approximations $\beta_\lambda$ and 
$\beta_{\Gamma,\lambda}$ (see, e.g., \cite[Appendix]{CF20b}). Hence, for the coupling term 
above we have the control
\[
  \int_\Gamma \beta_\lambda(v_{\delta,\lambda})\beta_{\Gamma,\lambda}(v_{\delta,\lambda})\,d\Gamma
  \geq
  \frac1{2\varrho_1}\int_{\Gamma}|\beta_\lambda(v_{\delta,\lambda})|^2 - C
\]
for some constant $C$. Then, integrating \eqref{8pier} over $(0,T)$ and applying the Young inequality, from (A1)--(A4) and the estimates \eqref{3rdf}, \eqref{6pier}, \eqref{1luca}, 
it is standard matter to deduce that
\begin{align}
&\lambda\int_\Omega\widehat\beta_\lambda(u_{\delta,\lambda}(T))\,dx
+ \lambda\int_\Gamma\widehat\beta_\lambda(v_{\delta,\lambda}(T))\,d\Gamma 
\nonumber \\
	&{} + \int_0^T \bigl| \beta_\lambda(u_{\delta,\lambda}(s))\bigr|_{H}^2 \,ds 
	+\int_0^T \bigl| \beta_\lambda(v_{\delta,\lambda}(s))\bigr|_{H_\Gamma}^2 \,ds 
	\le M_6
	\label{2luca}
\end{align}
for some positive constant $M_6$. Next, by comparing the terms in equation \eqref{vlw2a} we have that
\begin{equation*}
  \bigl| \Delta u_{\delta, \lambda}(t) \bigr|_{H} 
  \le \bigl| \mu_{\delta, \lambda}(t) \bigr| _{H} 
  + \lambda \bigl| \partial_t u_{\delta, \lambda}(t) 
  \bigl| _{H} + \bigl| \beta_\lambda \bigl( u_{\delta, \lambda}(t) \bigr) \bigr| _{H} 
  +\bigl| \pi \bigl( u_{\delta, \lambda}(t) \bigr) \bigr| _{H} + \bigl| f(t) \bigr|_{H}
\end{equation*}
for a.a.\ $t \in(0,T)$, whence
\begin{equation}
 \int_0^T  \bigl| \Delta u_{\delta, \lambda}(s) \bigr|_{H}^2 \,ds 
  \le M_6.
  \label{6th1}
\end{equation}

We proceed now by exploiting the idea of \cite{CFS21}.
Together with the trace theorems for the normal derivative, 
estimates \eqref{3rdf} and \eqref{6th1} yield that 
\begin{equation}
 \int_0^T  \bigl| \partial_{\boldsymbol{\nu}} u_{\delta,\lambda}(s)\bigr|_{Z_\Gamma^*}^2 \,ds 
  \le M_{6}.
  \label{3luca}
\end{equation} 
Analogously, recalling the estimate for 
$\delta^{1/2} \nabla_\Gamma v_{\delta,\lambda }$ in $L^\infty (0,T;H_\Gamma)$ in \eqref{3rdf}, 
by \eqref{vlw4a} the trace of $\delta^{1/2} u_{\delta,\lambda }$
is uniformly bounded in $L^2 (0,T;V_\Gamma)$. Therefore,
by virtue of the elliptic regularity (see, e.g., \cite[Theorem~3.2, p.\ 1.79]{BG87}) 
and again the trace theorems for the normal  derivative, we obtain that
\begin{equation}
  \delta \int_0^T \bigl| \partial_{\boldsymbol{\nu}} u_{\delta,\lambda}(s) \bigr|^2_{H_\Gamma} \,ds
  \le M_{7}.
  \label{6th2}
\end{equation}
Consequently, by comparing the terms in
equation \eqref{vlw6a} one deduces that
\begin{align*}
	& \bigl| -\delta \Delta_\Gamma v_{\delta, \lambda}(t) + 
	\beta_{\Gamma,\lambda} \bigl(v_{\delta,\lambda}(t) \bigr) \bigr|_{Z_\Gamma^*}
	\\ 
	& \le \bigl| \partial_{\boldsymbol{\nu}} u_{\delta, \lambda}(t) \bigr|_{Z_\Gamma^*} 
+ C \bigl( \bigl| w_{\delta, \lambda}(t) \bigr| _{H_\Gamma} 
	+ \lambda \bigl| \partial_t v_{\delta, \lambda}(t) 
	\bigl| _{H_\Gamma} + \bigl| \pi_\Gamma \bigl( v_{\delta, \lambda} (t) \bigr) \bigr| _{H_\Gamma} + \bigl| g(t) \bigr|_{H_\Gamma} \bigl)
\end{align*}
for a.a.\ $t \in(0,T)$, hence that
\begin{equation}
	\int_0^T  \bigl| -\delta \Delta_\Gamma v_{\delta, \lambda}(s) 
	+ \beta_{\Gamma,\lambda} \bigl(v_{\delta,\lambda}(s) \bigr) \bigr|_{Z_\Gamma^*}^2 \,ds 
	\le M_{7}.
	\label{7th}
\end{equation}
Since $\delta^{1/2}\Delta_\Gamma v_{\delta,\lambda}$ is bounded in $L^\infty(0,T;V_\Gamma^*)$
by the estimate \eqref{3rdf}, a direct comparison in \eqref{7th} yields also
\begin{equation}
	\delta \int_0^T  \bigl|  \Delta_\Gamma v_{\delta, \lambda}(s) \bigr|_{V_\Gamma^*}^2 \,ds 
	+ \int_0^T  \bigl| \beta_{\Gamma,\lambda} \bigl(v_{\delta,\lambda}(s) \bigr) \bigr|_{V_\Gamma^*}^2 \,ds 
	\le M_{7}.
	\label{7th2}
\end{equation}

\section{Proofs of main theorems}
\label{proofs}

We start by discussing the limiting procedure. 
The main issue concerns the passage to the limit as $\delta \to 0$. 
Indeed, it is known from \cite[Theorem~2.3]{CFW20} that letting 
$\lambda \to 0$ with weak and weak star convergences, we can prove Proposition~2.1. Moreover, 
the limit functions $u_\delta,\mu_\delta,\xi_\delta, v_\delta, w_\delta$, and $\eta_\delta$ satisfy 
\eqref{weak1}--\eqref{weak7} and same kind of uniform estimates obtained in the previous section, that is, 
the estimates
\begin{align}
	\sup_{t\in[0,T]}\bigl| u_{\delta}(t) \bigr|_{V}^2 
	+ \sup_{t\in[0,T]}
	\delta
	\bigl| \nabla_\Gamma v_{\delta}(t) \bigr|_{H_\Gamma}^2
	&\le M_3,
	\label{est1} \\
	\int_0^T \bigl| \partial _t u_{\delta}(s) \bigr|_{V_0^*}^2 \,ds
	+ 
	\int_0^T \bigl| \partial _t v_{\delta}(s) \bigr|_{V_{\Gamma,0}^*}^2 \,ds
	&\le M_3 
	\label{est2} \\
	\int_0^T \bigl| \mu_{\delta}(s) \bigr|^2_V \,ds + 
	\int_0^T \bigl| w_{\delta}(s) \bigr|^2_{V_\Gamma} \,ds 
	&\le 2 M_5, 
	\label{est3}\\
	\int_0^T\bigl| \xi_{\delta}(s) \bigr|_H ^2\,ds + 
	\int_0^t  \bigl| \Delta u_{\delta}(s) \bigr|_{H}^2 \,ds +
	\int_0^t \bigl| \partial_{\boldsymbol{\nu}} u_{\delta}(s) \bigr|^2_{Z_\Gamma^*} \,ds 
	&\le 3M_6
	\label{est4}\\
	\int_0^T  \bigl| -\delta \Delta_\Gamma v_{\delta}(s) 
	+ \eta_{\delta}(s) \bigr|_{Z_\Gamma^*}^2 \,ds 
	&\le M_{7},
	\label{est5}\\
	\delta \int_0^T \bigl|  \partial_{\boldsymbol{\nu}} u_{\delta}(s) \bigr|^2_{H_\Gamma} \,ds+
	\delta \int_0^T  \bigl|  \Delta_\Gamma v_{\delta}(s) \bigr|_{V_\Gamma^*}^2 \,ds 
	+ \int_0^T  \bigl| \eta_{\delta}(s)  \bigr|_{V_\Gamma^*}^2 \,ds 
	&\le 2M_7.
	\label{est6}
\end{align}
Moreover, we have that
\begin{equation}
	m \bigl(\partial_t u_{\delta}(t) \bigr) =0, \quad m_\Gamma \bigl(\partial_t v_{\delta,\lambda}(t) \bigr) =0
	\label{est9}
\end{equation}
for a.a.\ $t \in (0,T)$. As a remark, using \eqref{est2} and  
recalling the definition of norms in \eqref{norm}--\eqref{norm'}, 
we deduce similar uniform estimates 
for $\{ \partial _t u_{\delta} \}_{\delta \in (0,1]}$ in $L^2(0,T;V^*)$ and 
$\{ \partial _t v_{\delta} \}_{\delta \in (0,1]}$ in $L^2(0,T;V_\Gamma^*)$, respectively.

\subsection{Proof of Theorem~\ref{tak1}.} 
From \eqref{est1}--\eqref{est6} it follows that there exist a sextuplet $(u,\mu, \xi, v, w, \eta)$, with
\begin{align*}
	&u \in H^1(0,T;V^*) \cap L^\infty (0,T;V), \quad \Delta u \in L^2(0,T;H), \\
	&\mu \in L^2(0,T;V), \quad \xi \in L^2(0,T;H), \\
	&v \in H^1(0,T;V_\Gamma^*) \cap L^\infty (0,T;Z_\Gamma), \\
	&w \in L^2(0,T;V_\Gamma), \quad \eta \in L^2(0,T;Z_\Gamma^*),
\end{align*}
and a subsequence $\{\delta_k \}_{k \in \mathbb{N}}$ such that, as $k \to +\infty$, 
the convergences \eqref{cvpier1}--\eqref{cvpier10} hold.
Moreover, by virtue of the {A}ubin--{L}ions compactness results \cite{Sim87} and the 
compact embeddings
$V \mathop{\hookrightarrow} \mathop{\hookrightarrow} H$ 
and $Z_\Gamma \mathop{\hookrightarrow} \mathop{\hookrightarrow} H_\Gamma$, 
the following strong convergence properties
\begin{align}
	&u_{\delta_k} \to u \quad {\rm in~}C \takeshi{\bigl(} [0,T];H \takeshi{\bigr)}, 
	\label{conv11} \\
	&v_{\delta_k} \to v \quad {\rm in~}C \takeshi{\bigl(} [0,T];H_\Gamma \takeshi{\bigr)}
	\label{conv12}
\end{align}
hold as well.
The {L}ipschitz continuities of $\pi$ and $\pi_\Gamma$ give us then, as $k \to \infty$,
\begin{align}
	&\pi(u_{\delta_k}) \to \pi(u) \quad {\rm in~}C \takeshi{\bigl(}[0,T];H\takeshi{\bigr)}, 
	\label{conv13} \\
	&\pi_\Gamma(v_{\delta_k}) \to \pi_\Gamma(v) \quad {\rm in~}C\takeshi{\bigl(}[0,T];H_\Gamma\takeshi{\bigr)}.
	\label{conv14}
\end{align}
Therefore, taking the limit in \eqref{weak1} and \eqref{weak4} we can obtain 
the variational formulations \eqref{main1} and \eqref{main4}. The conditions 
\eqref{main7} and \eqref{main8} are also inferred from \eqref{weak6}--\eqref{weak7} 
on account of \eqref{conv11}--\eqref{conv12}.
Thanks to \eqref{cvpier1} and \eqref{cvpier6}, the boundary condition \eqref{main3} follows from \eqref{weak3} and the continuity of the linear trace operator 
$\gamma_0$ from $V$ to $Z_\Gamma$.

The first equation in \eqref{main2} is coming from the one in \eqref{weak2} 
owing to the convergences \eqref{cvpier2}, \eqref{cvpier4}, 
\eqref{cvpier5}, and \eqref{conv13}. 
The second condition in \eqref{main2} is proved by the 
demi-closedness of the maximal monotone operator induced by $\beta$, by virtue of
the strong convergence \eqref{conv11} and the weak convergence \eqref{cvpier5}. 
The variational formulation \eqref{main5} is also obtained from the first equation in \eqref{weak5}, 
due to the convergences \eqref{cvpier7}, \eqref{cvpier3}, 
\eqref{cvpier10}, and \eqref{conv14}. 
\smallskip

To conclude the proof of Theorem~\ref{tak1}, 
it remains to prove \eqref{main6}. 
To this aim, we multiply the equality in \eqref{weak2} by $u_{\delta_k}$ and integrate the resultant over $Q$ with respect to time and space variables. Using \eqref{weak3}, we have that
\begin{align}
	&\int_Q |\nabla u_{\delta_k} |^2 \,dxdt \,
	- \int_\Sigma \partial_{\boldsymbol{\nu}} u_{\delta_k} v_{\delta_k} \,d\Gamma dt 
	\nonumber \\
	&{}+
	\int_Q \xi_{\delta_k} u_{\delta_k} \,dxdt 
	=\int_Q \bigl( \mu_{\delta_k} -\pi(u_{\delta_k} )+f \bigr) u_{\delta_k} \,dxdt. 
	\label{last1}
\end{align}
On the other hand, multiplying the equality in \eqref{weak5} 
by $v_{\delta_k}$ and integrating then over $\Sigma$, we find out that
\begin{align}
	&\int_\Sigma \partial_{\boldsymbol{\nu}} u_{\delta_k} v_{\delta_k} \,d\Gamma dt 
	+
	{\delta_k} \int_\Sigma |\nabla_\Gamma v_{\delta_k} |^2 \,d\Gamma dt 
	\nonumber \\
	&{}+ 
	\int_\Sigma \eta_{\delta_k} v_{\delta_k} \,d\Gamma dt 
	=\int_\Sigma \bigl( w_{\delta_k} -\pi_\Gamma(v_{\delta_k} )+g \bigr) v_{\delta_k} \,d\Gamma dt. 
	\label{last2}
\end{align}
Summing \eqref{last1} and \eqref{last2}, using lower semicontinuity and weak-strong convergences, we obtain that 
\begin{align}
	& \limsup_{k \to +\infty} \int_\Sigma \eta_{\delta_k} v_{\delta_k} \,d\Gamma dt \nonumber \\
	& \le \limsup_{k \to +\infty} 
	\int_Q \bigl( \mu_{\delta_k} -\pi(u_{\delta_k} )+f \bigr) u_{\delta_k} \,dxdt + 
	\limsup_{k \to +\infty} 
	\int_\Sigma \bigl( w_{\delta_k} -\pi_\Gamma(v_{\delta_k} )+g \bigr) v_{\delta_k} \,d\Gamma dt 
	\nonumber \\
	& \quad {} - \liminf_{k \to +\infty} \int_Q |\nabla u_{\delta_k} |^2 \,dxdt  
	- \liminf_{k \to +\infty} \int_Q \xi_{\delta_k} u_{\delta_k} \,dxdt 
	- \liminf_{k \to +\infty}	{\delta_k} \int_\Sigma |\nabla_\Gamma v_{\delta_k} |^2 \,d\Gamma dt 
	\nonumber \\
	& \le \int_Q \bigl( \mu -\pi(u)+f \bigr) u \, dxdt + 
	\int_\Sigma \bigl( w -\pi_\Gamma(v)+g \bigr) v \,d\Gamma dt 
	\nonumber \\
	& \quad {} - \int_Q |\nabla u |^2 \, dxdt  
	- \int_Q \xi u \,dxdt  \,
	= \int_0^T \langle \eta, v \rangle _{Z_\Gamma^*,Z_\Gamma} \, dt
	\label{last3}
\end{align}  
and the last equality can be recovered combining the equation in \eqref{main2} tested by $u$ and \eqref{main5} with $z_\Gamma = v$ (cf.\ also \eqref{main3}).
Now, using the definition of subdifferential for $\beta_\Gamma $ in $L^2 (\Sigma)$, from the second inclusion in \eqref{weak5} we have that 
\begin{equation}
	\int_\Sigma \eta_{\delta_k} (\zeta_\Gamma - v_{\delta_k}) \,d\Gamma dt 
	+ \int_\Sigma \widehat{\beta}_\Gamma (v_{\delta_k})\,d\Gamma dt
	\le \int_\Sigma \widehat{\beta}_\Gamma (\zeta_\Gamma) \,d\Gamma dt 
	\label{liminf}
\end{equation} 
for all $\zeta_\Gamma \in L^2(0,T;H_\Gamma)$. 
If $\zeta_\Gamma \in L^2(0,T;V_\Gamma)$, then by virtue of the weak convergence \eqref{cvpier8}, 
the strong convergence \eqref{conv12}, 
the weak lower semicontinuity of $\widehat{\beta}_\Gamma$, and \eqref{last3}, we obtain
\begin{gather*}
	\lim _{k \to +\infty} \int_\Sigma \eta_{\delta_k} \zeta_\Gamma \,d\Gamma dt
	=\int_0^T \langle \eta, \zeta_\Gamma \rangle_{Z_\Gamma^*,Z_\Gamma} \,dt, 
	\\
	\liminf_{k \to +\infty} \left( - \int_\Sigma \eta_{\delta_k} v_{\delta_k} \,d\Gamma dt \right)
	= - \limsup_{k \to +\infty} \int_\Sigma \eta_{\delta_k} v_{\delta_k} \,d\Gamma dt 
	\ge - \int_0^T \langle \eta, v \rangle_{Z_\Gamma^*,Z_\Gamma}\, dt, 
	\\
	\liminf_{k \to +\infty} \int_\Sigma  \widehat{\beta}_\Gamma (v_{\delta_k})\, d\Gamma dt 
	\ge 
	\int_\Sigma  \widehat{\beta}_\Gamma (v) \,d\Gamma dt.
\end{gather*}
Therefore, taking the infimum limit in \eqref{liminf}, we deduce that
\begin{equation}
	\int_0^T \langle \eta, \zeta_\Gamma-v \rangle_{Z_\Gamma^*,Z_\Gamma } \,dt 
	\le \int_\Sigma \widehat{\beta}_\Gamma (\zeta_\Gamma) \,d\Gamma dt -  
	\int_\Sigma \widehat{\beta}_\Gamma (v) \,d\Gamma dt 
	\label{sub}
\end{equation}
for all $\zeta_\Gamma \in L^2(0,T;V_\Gamma)$. 
As $\eta \in L^2 (0,T;Z_\Gamma^*)$, by a density argument we can prove that \eqref{sub} holds also for all $\zeta_\Gamma \in L^2(0,T;Z_\Gamma)$. Indeed,
for a given arbitrary $\zeta_\Gamma \in L^2(0,T;Z_\Gamma)$ and $\varepsilon>0$, 
we can take the approximations 
$\{ \zeta_{\Gamma,\varepsilon} \}_{\varepsilon>0} \subset L^2(0,T;V_\Gamma)$ 
defined as the solutions to
\begin{equation*}
	\zeta_{\Gamma, \varepsilon} - \varepsilon \Delta_\Gamma \zeta_{\Gamma,\varepsilon} =\zeta_\Gamma 
	\quad {\rm a.e.\ on~} \Sigma.  
\end{equation*}
In fact, thanks to \cite[Lemma~A.1]{CF16} we have that
\begin{gather*}
	\zeta_{\Gamma,\varepsilon} \to \zeta_\Gamma 
	\quad {\rm in~} L^2(0,T;Z_\Gamma) \quad {\rm as~} \varepsilon \to 0, 
	\\
	\widehat{\beta}_\Gamma (\zeta_{\Gamma,\varepsilon}) \le \widehat{\beta}_\Gamma (\zeta_\Gamma) \quad 
	{\rm a.e.\ on~} \Sigma, \hbox{ for~all~} \varepsilon >0.
\end{gather*}
Thus, replacing $\zeta_\Gamma$ by $\zeta_{\Gamma,\varepsilon}$ in \eqref{sub} 
and letting $\varepsilon \to 0$, we 
easily obtain the validity of \eqref{sub} for all $\zeta_\Gamma \in L^2(0,T;Z_\Gamma)$, which is an equivalent formulation of \eqref{main6}. \hfill $\Box$

\subsection{Proof of Theorem~\ref{tak2}.}
Next, we prove the continuous dependence result stated in Theorem~\ref{tak2}.
Assume that $f^{(1)}, f^{(2)}$, $g^{(1)}, g^{(2)}$ satisfy {\rm (A3)}, 
$u_{0}^{(1)}, u_{0}^{(2)}, v_{0}^{(1)},v_{0}^{(2)}$ satisfy {\rm (A4)} and 
\begin{equation}
	m\bigl(u_{0}^{(1)}\bigr)=m\bigl(u_{0}^{(2)}\bigr), \quad m_\Gamma\bigl(v_{0}^{(1)}\bigr)=m_\Gamma\bigl(v_{0}^{(2)}\bigr).
	\label{ad}
\end{equation}
For these data  let 
$(u^{(i)},\mu^{(i)}, \xi^{(i)}, v^{(i)}, w^{(i)}, \eta^{(i)})$, $i=1,2$, 
be respective solutions obtained by 
Theorem~\ref{tak1} 
Put $\bar{u}:=u^{(1)}-u^{(2)}$ and analogously use the same notation for the differences of functions.
Taking the difference of \eqref{main1}, \eqref{main2}, \eqref{main4}, and \eqref{main5}, we have  
\begin{gather} 
	\langle \partial_t \bar{u} , z \rangle_{V^*,V}+ \int_\Omega \nabla \bar{\mu}  \cdot \nabla z \,dx=0,
	\label{con1}
	\\
	(\bar{\mu}, z)_H = (\nabla \bar{u}, \nabla z)_{H}
	-\langle \partial_{\boldsymbol{\nu}} \bar{u}, z_{|_\Gamma} \rangle_{Z_\Gamma^*,Z_\Gamma} 
	+ (\bar{\xi},z)_H + \bigl( \pi (u_1)-\pi(u_2) -\bar{f},z \bigr)_{\!H}
	\label{con2}
\end{gather}
for all $z \in V$ and a.e.\ in $(0,T)$,
\begin{gather} 
	\langle \partial_t \bar{v}, z_\Gamma  \rangle_{V_\Gamma^*,V_\Gamma}
	+ \int_{\Gamma} \nabla_\Gamma \bar{w} \cdot \nabla_\Gamma z_\Gamma \,d\Gamma = 0,
	\label{con3}
\end{gather}
for all $z_\Gamma \in V_\Gamma$ and a.e.\ in $(0,T)$,
\begin{gather}
	(\bar{w}, z_\Gamma)_{H_\Gamma} = \langle \partial_{\boldsymbol{\nu}} \bar{u}+\bar{\eta}, z_\Gamma \rangle_{Z_\Gamma^*,Z_\Gamma}   
	+ \bigl( \pi_\Gamma(v_1)-\pi(v_2)-\bar{g},z_\Gamma \bigr)_{\!H_\Gamma},
	\label{con4}
\end{gather}
for all $z_\Gamma \in Z_\Gamma$ 
and a.e.\ in $(0,T)$. Moreover, using \eqref{ad} we have 
\begin{equation*}
	m \bigl(\bar{u}(t) \bigr)=0, \quad m_\Gamma \bigl( \bar{v}(t)\bigr)=0
\end{equation*}
for all $t \in [0,T]$. 
Take $z=F^{-1}\bar{u}$ in \eqref{con1}, $z=\bar{u}$ in \eqref{con2}, 
$z_\Gamma =F_\Gamma^{-1} \bar{v}$ in \eqref{con3}, and $z_\Gamma=\bar{v}$ in \eqref{con4}, respectively. 
Then, combining the obtained equalities and integrating over $(0,t)$, we deduce that
\begin{align*}
	& \frac{1}{2} \bigl| \bar{u}(t) \bigr|_*^2 
	+ 
	\frac{1}{2} \bigl| \bar{v}(t) \bigr|_{\Gamma,*}^2
	+ \int_0^t \bigl| \bar{u}(s) \bigr|_{V_0}^2 \,ds 
	+ \int_0^t \bigl( \bar{\xi}(s),\bar{u}(s)\bigr) _H \,ds 
	+ \int_0^t \bigl \langle \bar{\eta}(s),\bar{v}(s)\bigr \rangle_{Z_\Gamma^*,Z_\Gamma} \,ds 
	\nonumber \\
	& = \frac{1}{2} | \bar{u}_0 |_*^2 
	+ 
	\frac{1}{2} | \bar{v}_0 \bigr|_{\Gamma,*}^2
	+
	\int_0^t \bigl( \bar{f}(s)+\pi \bigl(u^{(2)}(s) \bigr)-\pi \bigl(u^{(1)}(s) \bigr) , \bar{u}(s) \bigr)_{\!H} \,ds 
	\nonumber \\
	& \quad {} + \int_0^t \bigl( \bar{g}(s)+\pi_\Gamma \bigl(v^{(2)}(s) \bigr)
	-\pi_\Gamma \bigl(v^{(1)}(s) \bigr) , \bar{v}(s) \bigr)_{\!H_\Gamma}\, ds 
\end{align*}
for all $t \in [0,T]$. 
Now, we invoke the monotonicity of the maximal monotone operators induced by $\beta$ and $\beta_\Gamma$ (cf.~Remark~\ref{rem1}) to see that the last two terms on the left-hand side are nonnegative. We also use the following estimate
\begin{equation*}
	\bigl| \bar{v}(s) \bigr|_{Z_\Gamma}^2 \le C_{\rm tr}^2 \bigl| \bar{u}(s) \bigr|_V^2 \le C_{\rm tr}^2 C_{\rm P} \bigl| \bar{u}(s) \bigr|_{V_0}^2,
\end{equation*}
which comes from \eqref{tr0} and \eqref{pier4}. Then, on account of 
 the {L}ipschitz continuity of $\pi$ and $\pi_\Gamma$, by applying twice
the {E}hrling lemma we can conclude that for all $\varepsilon>0$ there is a constant 
$C_\varepsilon >0 $ such that
\begin{align*}
	& \bigl| \bar{u}(t) \bigr|_*^2 
	+ 
	\bigl| \bar{v}(t) \bigr|_{\Gamma,*}^2
	+ \int_0^t \bigl| \bar{u}(s) \bigr|_{V_0}^2 \,ds 
	+ \frac{1}{ C_{\rm tr}^2 C_{\rm P}}
	\int_0^t \bigl| \bar{v}(s) \bigr|_{Z_\Gamma}^2\, ds 
	\nonumber \\
	& \le | \bar{u}_0 |_*^2 
	+ 
	| \bar{v}_0 \bigr|_{\Gamma,*}^2
	+
	\int_0^t \bigl| \bar{f}(s) \bigr|_H^2 \,ds
	+ \varepsilon \int_0^t \bigl| \bar{u}(s) \bigr|_{V_0}^2 \,ds +
	 C_\varepsilon \int_0^t \bigl| \bar{u}(s) \bigr|_{*}^2 \,ds
	\nonumber \\
	& \quad {} + \int_0^t \bigl| \bar{g}(s) \bigr|_{H_\Gamma}^2 \,ds
	+ \varepsilon \int_0^t \bigl| \bar{v}(s) \bigr|_{Z_\Gamma}^2 \,ds +
	 C_\varepsilon \int_0^t \bigl| \bar{v}(s) \bigr|_{\Gamma,*}^2 \,ds
\end{align*}
for all $t \in [0,T]$ . Thus, choosing $\varepsilon>0$ sufficiently 
small and applying the Gronwall lemma, by the {P}oincar\'e--Wirtiger inequality \eqref{pier4} we complete the proof of Theorem~\ref{tak2}. \hfill $\Box$

\subsection{Proof of Theorem~\ref{tak3}.}
We point out that the further assumption \eqref{ip_extra} on the graphs
yields additional estimates on the solutions. Indeed, 
since assumption \eqref{ip_extra} induces the same 
inequalities on the respective Yosida approximations
(details are given in \cite[Appendix]{CF20}) and, in particular, \eqref{ip_extra}
implies that 
\[
  \frac1{2 M^2}\int_\Gamma 
  \takeshi{\bigl|} 
  \beta_{\Gamma,\lambda}
  \takeshi{\bigl(} 
  v_{\delta,\lambda}
  \takeshi{(t)}
  \takeshi{\bigr)} \takeshi{\bigr|}^2 \, d\Gamma
  \leq
  \int_{\Gamma} \bigl(
  \takeshi{\bigl|} 
  \beta_\lambda 
  \takeshi{\bigl(} 
  v_{\delta,\lambda} \takeshi{(t)}  
  \takeshi{\bigr)} 
  \takeshi{\bigr|}^2 + M^2 \bigr)
	\, \takeshi{d\Gamma}
\]
\takeshi{for a.a.\ $t \in (0,T)$,} 
the estimate \eqref{2luca} entails that
\[
  \takeshi{\bigl|} 
  \beta_\lambda(u_{\delta,\lambda})
  \takeshi{\bigr|}_{L^2(0,T; H)}+
  \takeshi{\bigl|}
  \beta_{\Gamma,\lambda}(v_{\delta,\lambda})
  \takeshi{\bigr|}_{L^2(0,T; H_\Gamma)}\leq C
  \]
for some positive constant $C$. Hence, recalling the
estimates \eqref{3rdf}, \eqref{1luca} and \eqref{3luca}, by comparison of terms in \eqref{vlw6a} 
we find out that
\begin{equation}\label{est14}
  \takeshi{\bigl|} 
  \partial_{\boldsymbol{\nu}} u_{\delta,\lambda} - \delta\Delta_\Gamma v_{\delta,\lambda}\takeshi{\bigr|}_{L^2(0,T; H_\Gamma)}
  +\delta \takeshi{\bigl|} 
  \Delta_\Gamma v_{\delta,\lambda}
  \takeshi{\bigr|}_{L^2(0,T; Z_\Gamma^*)} \leq C
\end{equation}
and consequently, by elliptic regularity, 
\begin{equation}\label{est14-bis}
  \takeshi{\bigl|} 
  \delta v_{\delta,\lambda} \takeshi{\bigr|}_{L^2(0,T; H^{3/2}(\Gamma))} \leq C . 
\end{equation}
Then, we can take the limit as $\lambda\to0$ and infer that 
\begin{align} 
&\int_0^T  \bigl| \eta_{\delta}(s) \bigr|_{H_\Gamma}^2 \,ds 
+ \int_0^T \bigl|  \partial_{\boldsymbol{\nu}} u_{\delta}(s) - \delta \Delta_\Gamma v_{\delta}(s) \bigr|^2_{H_\Gamma} \,ds
+	\int_0^T  \bigl| \delta v_{\delta}(s) 
	\bigr|_{H^{3/2}(\Gamma)}^2 \,ds 
	\le C 
	\label{est5-bis}
	\end{align}
in addition to \eqref{est1}--\eqref{est6}. Thus, in view of \eqref{cvpier1}--\eqref{cvpier10}, when passing to the limit on a subsequence $\delta_k$ we also deduce 
\eqref{cvpier11}--\eqref{cvpier13} and the boundary equation \eqref{pier7} at the limit. 
At this point, as $u\in L^2(0,T; V)$, $\Delta u \in  L^2(0,T; H)$  and  $\partial_{\boldsymbol{\nu}} u \in  L^2(0,T; H_\Gamma)$,
by elliptic regularity (see \cite[Thm.~3.2]{BG87}) it follows that 
\[  
  u  \in L^2 \takeshi{\bigl(} 0,T; H^{3/2}(\Omega) \takeshi{\bigr)},
\]
whence, from \eqref{main3} and the trace theory, 
\[
  v \in L^2(0,T; V_\Gamma).
\]
Eventually, the pointwise inclusion $\xi_\Gamma\in\beta_\Gamma(u_\Gamma)$
a.e.\ on $\Sigma$ is ensured in this framework, as explained in 
Remark~\ref{rem1}. This ends the proof of Theorem~\ref{tak3}.\hfill $\Box$

\subsection{Proof of Theorem~\ref{tak4}.} 
For $\delta \in (0,1]$ let $(u_\delta,
\mu_\delta, \xi_\delta, v_\delta, w_\delta, \eta_\delta)$ be the sextuplet, solution of the 
problem~\eqref{weak1}--\eqref{weak7}, obtained in the passage to the limit as $
\lambda\to 0$ and let $(u,\mu, \xi, v, w, \eta)$ denote the solution of the problem~\eqref{main1}--\eqref{main8} arising from the above proof of Theorem~\ref{tak3} (cf.\ Theorem~\ref{tak1} as well).

Now, we argue similarly as in the proof of Theorem~\ref{tak2} and use the notations 
$\bar{u}_\delta:=u_\delta - u$, $\bar{\mu}_\delta:=\mu_\delta - \mu$, 
$\bar{\xi}_\delta:=\xi_\delta - \xi$,
$\bar{v}_\delta:=v_\delta - v$, $\bar{w}_\delta:=w_\delta - w$, $\bar{\eta}_\delta:=\eta_\delta - \eta$. Here, in place of \eqref{con1}--\eqref{con4} we have 
the equalities 
\begin{gather} 
	\langle \partial_t \bar{u}_\delta , z \rangle_{V^*,V}+ \int_\Omega \nabla \bar{\mu}_\delta  \cdot \nabla z \,dx=0,
	\label{conpi1}
	\\
	(\bar{\mu}_\delta, z)_H = (\nabla \bar{u}_\delta, \nabla z)_{H}
	- ( \partial_{\boldsymbol{\nu}} \bar{u}_\delta, z_{|_\Gamma} )_{H_\Gamma}
	+ (\bar{\xi}_\delta,z)_{H} + \bigl( \pi (u_\delta)-\pi(u) ,z \bigr)_{\!H}
	\label{conpi2}
\end{gather}
for all $z \in V$ and a.e.\ in $(0,T)$;
\begin{gather} 
	\langle \partial_t \bar{v}_\delta, z_\Gamma  \rangle_{V_\Gamma^*,V_\Gamma}
	+ \int_{\Gamma} \nabla_\Gamma \bar{w}_\delta \cdot \nabla_\Gamma z_\Gamma \,d\Gamma = 0,
	\label{conpi3}
\\	(\bar{w}_\delta, z_\Gamma)_{H_\Gamma} = 
\delta \int_{\Gamma} \nabla_\Gamma v_\delta \cdot \nabla_\Gamma z_\Gamma \,d\Gamma
	+(\partial_{\boldsymbol{\nu}} \bar{u}_\delta + \bar{\eta}_\delta, z_\Gamma )_{H_\Gamma}  
	+ \bigl( \pi_\Gamma(v_\delta)-\pi(v),z_\Gamma \bigr)_{\!H_\Gamma}
	\label{conpi4}
\end{gather}
for all $z_\Gamma \in V_\Gamma$ and a.e.\ in $(0,T)$. As
\begin{equation*}
	m \bigl(\bar{u}_\delta (s) \bigr)=0, \quad m_\Gamma \bigl(  \bar{u}_\delta(s)\bigr)=0
\end{equation*}
for all $s \in [0,T]$, 
we can take $z=F^{-1}\bar{u}_\delta\takeshi{(s)}$ in \eqref{conpi1}, 
$z=- \bar{u}_\delta\takeshi{(s)}$ in \eqref{conpi2}, 
and add them with a cancellation; then,
we choose $z_\Gamma =F_\Gamma^{-1} \bar{v}_\delta\takeshi{(s)}$ in \eqref{conpi3}, 
and $z_\Gamma= - \bar{v}_\delta\takeshi{(s)}$ in \eqref{conpi4}, and add the two resultants with another cancellation. Finally, we can take the sum and integrate over $(0,t)$, obtaining  
\begin{align*}
	& \frac{1}{2} \bigl| \bar{u}_\delta(t) \bigr|_*^2 
	+ \frac{1}{2} \bigl| \bar{v}_\delta(t) \bigr|_{\Gamma,*}^2
	+ \int_0^t \bigl| \bar{u}_\delta(s) \bigr|_{V_0}^2 \,ds 
	+ \delta \int_0^t \!\!\int_{\Gamma} \takeshi{\bigl|} 
	\nabla_\Gamma v_\delta (s) \takeshi{\bigr|}^2 \, d\Gamma ds 	
	\\
	&+ \int_0^t \bigl( \bar{\xi}_\gamma(s),\bar{u}_\gamma(s)\bigr)_{\!H} \,ds 
	+ \int_0^t \bigl( \bar{\eta}_\gamma(s),\bar{v}(s) \bigr)_{\!H_\Gamma} \,ds 
	\\
	& = \delta \int_0^t \!\!\int_{\Gamma} \nabla_\Gamma v_\delta (s) \cdot
	\nabla_\Gamma v (s)	\, d\Gamma ds	+
	\int_0^t \bigl( \pi \bigl(u_\delta(s) \bigr)-\pi \bigl(u(s) \bigr) , \bar{u}_\delta(s) \bigr)_{\!H} \,ds 
	\\
	& \quad {} + \int_0^t \bigl( \pi_\Gamma \bigl(v_\delta (s) \bigr)
	-\pi_\Gamma \bigl(v(s) \bigr) , \bar{v}_\delta(s) \bigr)_{\!H_\Gamma}\, ds 
\end{align*}
for all $t \in [0,T]$. Next, we observe that
$$ \int_0^t \bigl( \bar{\xi}_\gamma(s),\bar{u}_\gamma(s)\bigr)_{\!H} \,ds \geq 0 ,\quad
 \int_0^t \bigl( \bar{\eta}_\gamma(s),\bar{v}(s) \bigr)_{\!H_\Gamma} \,ds \geq 0$$
due to the monotonicity of $\beta $ and $\beta_\Gamma$;
$$ \delta \int_0^t \!\!\int_{\Gamma} \nabla_\Gamma v_\delta (s) \cdot
	\nabla_\Gamma v (s)	\, d\Gamma ds 
	\leq \frac\delta2 \int_0^t \!\!\int_{\Gamma} 
	\takeshi{\bigl|} 
	\nabla_\Gamma v_\delta (s)
	\takeshi{\bigr|}^2 \, d\Gamma ds
	+ \frac\delta2 \int_0^t \!\!\int_{\Gamma} 
	\takeshi{\bigl|}\nabla_\Gamma v (s) \takeshi{\bigr|}^2 \, d\Gamma ds
$$
by the Young inequality; moreover, we can treat the terms containing the differences 
\takeshi{$\pi (u_\delta(s))-\pi(u(s))$ and 
$\pi_\Gamma (v_\delta (s))-\pi_\Gamma (v(s))$} exactly in the same way as in the proof of Theorem~\ref{tak2}, using Lipschitz continuity and the Ehrling lemma. Then, with the help of the Gronwall lemma and the {P}oincar\'e--Wirtiger inequality~\eqref{pier4} we arrive at 
\begin{align*}
	& \takeshi{|} \bar{u}_\delta \takeshi{|}_{L^\infty(0,T;V^*)}^2 
	+ \takeshi{|} \bar{v}_\delta \takeshi{|}_{L^\infty(0,T;V_\Gamma^*)}^2 	
	+ \takeshi{|}  \bar{u}_\delta \takeshi{|}_{L^2(0,T;V)}^2 \\
	&+ \takeshi{|} \bar{v}_\delta \takeshi{|}_{L^2(0,T;Z_\Gamma)}^2 			
	+ \delta \int_0^T \!\!\int_{\Gamma} \takeshi{\bigl|} 
	\nabla_\Gamma v_\delta (t) \takeshi{\bigr|}^2 \, d\Gamma dt
\leq C\delta \int_0^T \!\!\int_{\Gamma} \takeshi{\bigl|} 
\nabla_\Gamma v (t)\takeshi{ \bigr|}^2 \, d\Gamma dt
\end{align*} 
for some positive constant $C$ depending only on data. Then, as $v$ belongs to $L^2(0,T;V_\Gamma)$, it is straightforward to deduce both the error estimate \eqref{pier8}
and the additional convergence \eqref{pier9}, which is a consequence of the boundedness of
$ \int_0^T \!\!\int_{\Gamma} |\nabla_\Gamma v_\delta (t)|^2 d\Gamma dt$ independent of $\delta$ and the strong convergence ${v}_\delta \to v$ in $L^2(0,T;Z_\Gamma)$.\hfill$\Box$

\section*{Acknowledgments}
This research received a support from the Italian Ministry of Education, 
University and Research~(MIUR): Dipartimenti di Eccellenza Program (2018--2022) 
-- Dept.~of Mathematics ``F.~Casorati'', University of Pavia.
TF acknowledges the support from the JSPS KAKENHI Grant-in-Aid for Scientific Research(C), 
Japan, Grant Number 17K05321 and from the Grant Program of The Sumitomo Foundation, Grant 
Number 190367. PC and LS gratefully acknowledge some support 
from the MIUR-PRIN Grant 2020F3NCPX ``Mathematics for industry 4.0 (Math4I4)'' 
and underline their affiliation
to the GNAMPA (Gruppo Nazionale per l'Analisi Matematica, 
la Probabilit\`a e le loro Applicazioni) of INdAM (Isti\-tuto 
Nazionale di Alta Matematica). Moreover, PC aims to point out his collaboration,
as Research Associate, to the IMATI -- C.N.R. Pavia, Italy.

\end{document}